\documentclass[12pt]{amsart}

\usepackage{tcolorbox}
\usepackage{xcolor}
\usepackage{extpfeil}

\usepackage{pdfsync}
\usepackage{geometry}
\usepackage{lscape}
\usepackage{graphicx}
\usepackage{epstopdf}    
\usepackage[all]{xy}
\usepackage{amsmath,amsthm,amssymb,enumerate}
\usepackage{latexsym}
\usepackage{amscd}
\usepackage{amssymb}
\usepackage{enumitem}
\usepackage{mathtools}

\numberwithin{equation}{section}

\theoremstyle{plain}  % These will have emph text
\newtheorem{thm}[equation]{Theorem}
\newtheorem{prop}[equation]{Proposition}

\theoremstyle{definition}  % These will have plain text

\newtheorem{remark}[equation]{Remark}

\newcommand{\ul}{\underline}

\newcommand{\ra}{\rightarrow}

\newcommand{\lra}{\longrightarrow}

\newcommand{\llra}[1]{\stackrel{#1}{\lra}}

\newcommand{\KO}[1]{ \ul{KO}_{\;#1}}
\newcommand{\HKO}[1]{ H_*(\ul{KO}_{\;#1})}

\newcommand{\KU}[1]{ \ul{KU}_{\;#1}}
\newcommand{\KR}[1]{ \ul{KR(1)}_{\;#1}}
\newcommand{\HKU}[1]{H_*( \ul{KU}_{\;#1})}

\newcommand{\HKR}[1]{ H_*(\ul{KR(1)}_{\;#1})}

\newcommand{\K}[1]{ \ul{K(1)}_{\;#1}}
\newcommand{\HK}[1]{ H_*(\ul{K(1)}_{\;#1})}

\newcommand{\Z}{\mathbb Z}

\newcommand{\Zq}{\Z/(2)}

\DeclareMathOperator{\Tor}{Tor}

\thanks{
}

\begin{document}

\title{The Omega spectrum for mod 2 $KO$-theory }

\author{W. Stephen Wilson}

\address{Department of Mathematics, The Johns Hopkins University, Baltimore, MD 21218}

\email{wwilson3@jhu.edu}

\begin{abstract}
	The 8-periodic theory that comes from the KO-theory of the mod 2 Moore
	space is the same as the real first Morava K-theory obtained from the
	homotopy fixed points of the Z/(2) action on the first Morava K-theory.
	The first Morava K-theory, K(1), is just mod 2 KU-theory.
	We compute the homology Hopf algebras for the spaces in this Omega spectrum.
There are a lot of maps
into and out of these spaces and the spaces for KO-theory, KU-theory
and the first Morava K-theory.  For every one of these 98 maps (counting
suspensions) there is a spectral sequence.  
We describe
all 98 maps and spectral sequences.  48 of these maps involve
our new spaces and 56 of the spectral sequences do.
In addition, the maps on homotopy are all written down.
\end{abstract}

\thanks{Thanks to Nitu Kitchloo, Don Davis, and Bill Singer, for some help with this.}

\date{\today}

\maketitle

\section{Introduction}

We have stable maps $2 \colon S^0 \ra S^0$ and $\eta \colon S^1 \ra S^0$ and we
get a stable diagram with $M$ the mod 2 Moore space and $N$ and $NM$ the
appropriate cofibres
\[
	\xymatrix{
		S^1 \ar[r]^2 \ar[d]^\eta & S^1 \ar[d]^\eta \ar[r]  & \Sigma^1 M \ar[d]^\eta \\
		S^0 \ar[r]^2 \ar[d] & S^0 \ar[r] \ar[d]  & M  \ar[d]\\
		N \ar[r]^2 & N \ar[r] & NM
		}
	\]

	If we smash this diagram with connective $K$-theory, $bo$, and then only look
	at the low dimensional spaces in the Omega spectrum where we get periodicity,
	we get the diagram of fibrations

\begin{equation}
\xymatrix{
	\KO{i+1}\ar[d]^\eta \ar[r]^2 & 
\KO{i+1}\ar[d]^\eta \ar[r]^\rho & \KR{i+1} \ar[d]^\eta \ar[r]^\delta & \ar[d]^\eta \KO{i+2} \\
\KO{i}\ar[d]^\rho \ar[r]^2 & \KO{i}\ar[d]^\rho \ar[r]^\rho & 
\KR{i} \ar[d]^\rho \ar[r]^\delta & \KO{i+1}  \ar[d]^\rho     \\
	\KU{i} \ar[d]^\delta  \ar[r]^2 & 
\KU{i} \ar[d]^\delta  \ar[r]^\rho  & 
\K{i} \ar[d]^\delta  \ar[r]^\delta & \ar[d]^\delta  \KU{i+1}  \\
	\KO{i+2} \ar[r]^2 & 
\KO{i+2} \ar[r]^\rho & \KR{i+2}  \ar[r]^\delta &  \KO{i+3} \\
	}
\end{equation}

The $\KU{i}$ are the usual 2-periodic
spaces for complex $K$-theory and the $\KO{i}$ the 8-periodic spaces for
real
$K$-theory.  The $\K{i}$ are 2-periodic and they are just the mod 2
$KU$-theory, or, the first Morava $K$-theory.
The spaces of interest are the $\KR{i}$, which are simultaneously the real version
of the first Morava $K$-theory (see \cite[Theorem 3.32]{HK}) and the mod 2 $KO$-theory.

Our interest is in computing the Hopf algebra $\HKR{i}$.
We work with $\Zq$ coefficients in homology.
Our notation is that $P$ is a polynomial algebra, $E$ is an exterior
algebra, $TP_4(x)$ is $P(x)/(x^4)$.  The Frobenius, $F$, is just
the map that takes $x$ to $x^2$.  The Vershiebung is the dual of the
Frobenius and gives us the coproduct structure on our Hopf algebras.
Our notation is such that the subscript of an element denotes the
degree it resides in.

Our main theorem is easy to state.

\begin{thm} The homology of the connected component of $\KR{i}$ is as follows.
	If the Vershiebung isn't described, it is zero.
The index $k$ runs over all $k > 0$.
\begin{alignat*}{4}
&i=0\hspace{.5in} 
&
E(x_k) &\otimes_k P(y_{4k+2})
&
\hspace{.25in}V(x_{2k}) &= x_k
& \\
&i=1\hspace{.5in}
&
P(x_{2k+1})&\otimes_k P(y_{4k+2}) 
&
\hspace{.55in}V(y_{4k+2}) &= x_{2k+1}
&\\
&i=2\hspace{.5in}
&
P(x_{8k+2}) &\otimes_k P(y_{4k+3})
&
&
& \\
&i=3
&
E(x_{8k+3}) &\otimes_k P(y_{8k+4})
&
&
& \\
&i=4
&
E(x_{4k}) &\otimes_k E(y_{8k+5})
&
\hspace{.25in}V(x_{8k}) &= x_{4k}
& \\
&i=5
&
E(x_{4k+1}) &\otimes_k E(y_{2k})
&
\hspace{.25in}V(y_{4k}) &= y_{2k}
\qquad V(y_{8k+2}) &=  x_{4k+1}
\\
&i=6
&
\otimes_k &TP_4(x_{k})
&
\hspace{.25in}V(x_{2k}) &= x_{k}
& \\
&i=7
&
E(x_{2k})
&\otimes_k 
P(y_{2k+1}) 
&
\hspace{.25in}V(x_{4k}) &= x_{2k}
&
\end{alignat*}
%\[
%	\begin{array}{lcll}
%i=0\hspace{.5in} &
%E(x_k) \otimes P(y_{4k+2})
%&
%\hspace{.25in}V(x_{2k}) = x_k
%& \\
%i=1\hspace{.5in}
%&
%P(x_{2k+1})\otimes P(y_{4k+2}) 
%&
%\hspace{.25in}V(y_{4k+2}) = x_{2k+1}
%&
%\\
%i=2
%&
%P(x_{8k+2}) \otimes P(y_{4k+3})
%&
%&
%\\
%i=3
%&
%E(x_{8k+3}) \otimes P(y_{8k+4})
%&
%& \\
%i=4
%&
%E(x_{4k}) \otimes E(y_{8k+5})
%&
%\hspace{.25in}V(x_{8k}) = x_{4k}
%& \\
%i=5
%&
%E(x_{4k+1}) \otimes E(y_{2k})
%&
%\hspace{.25in}V(y_{4k}) = y_{2k}
%&
%V(y_{8k+2}) =  x_{4k+1}
%\\
%i=6
%&
%TP_4(x_{k})
%&
%\hspace{.25in}V(x_{2k}) = x_{k}
%& \\
%i=7
%&
%E(x_{2k})
%\otimes 
%P(y_{2k+1}) 
%&
%\hspace{.25in}V(x_{4k}) = x_{2k}
%&
%\\
%\end{array}
%\]
\end{thm}

%\begin{thm} The homology of the connected component of $\KR{i}$ is as follows.
%	If the Vershiebung isn't described, it is zero.
%	\begin{itemize}
%		\item[i=0]
%\qquad
%$E(x_k) \otimes P(y_{4k+2})
%\qquad
%V(x_{2k}) = x_k
%$
%\item[i=1]
%\qquad
%$P(x_{2k+1})\otimes P(y_{4k+2}) 
%\qquad
%V(y_{4k+2}) = x_{2k+1}
%$
%\item[i=2]
%\qquad
%$P(x_{8k+2}) \otimes P(y_{4k+3})
%$
%\item[i=3]
%\qquad
%$E(x_{8k+3}) \otimes P(y_{8k+4})
%$
%\item[i=4]
%\qquad
%$E(x_{4k}) \otimes E(y_{8k+5})
%\qquad
%V(x_{8k}) = x_{4k}
%$
%\item[i=5]
%\qquad
%$
%E(x_{4k+1}) \otimes E(y_{2k})
%\qquad
%V(y_{4k}) = y_{2k}
%\qquad
%V(y_{8k+2}) =  x_{4k+1}
%$
%\item[i=6]
%\qquad
%$TP_4(x_{k})
%\qquad
%\qquad
%\qquad
%V(x_{2k}) = x_{k}
%$
%\item[i=7]
%\qquad
%$
%E(x_{2k})
%\otimes 
%P(y_{2k+1}) 
%\qquad
%V(x_{4k}) = x_{2k}
%$
%	\end{itemize}
%\end{thm}

\begin{remark}
We began this research trying to give meaningful names to all of
the algebra generators.  Eventually, it became clear that it was
easier to compute just using the degrees of the generators.
We do know good names for all of the generators of 
$\HKO{*}$, $\HKU{*}$, and $\HK{*}$, and we are able to relate
our poorly named generators to generators we are more familiar with,
thus solving the naming problem after the fact.
In order to be explicit about these results, we have to write down the known
homologies first.  We will put that off until the next section.
We can give the non-explicit answer here.
\end{remark}

\begin{thm}
\label{exactness}
The maps of the connected components $\KO{i} \xrightarrow{\rho}
\KR{i} \xrightarrow{\delta} \KO{i+1}$ give rise to maps on homology
\[
\HKO{i} \xrightarrow{\rho_*}
\HKR{i} \xrightarrow{\delta_*} \HKO{i+1}
\]
that are exact in the category of Hopf algebras at the middle term.
For $i=1,2,5$ and $6$, this is a short exact sequence of Hopf algebras.
In the case $i=0$ we have a long exact sequence:
%\begin{align*}
%\HKO{1}\xrightarrow[\hspace{.7cm}]{\eta} \HKO{0} \xrightarrow[\hspace{.7cm}]{2} \HKO{0} 
%\xrightarrow[\hspace{.7cm}]{\rho} \HKR{0} \xrightarrow[\hspace{.7cm}]{\delta} \HKO{1}
%\xrightarrow[\hspace{.7cm}]{\rho} \HKU{1}
%\end{align*}
\[
 \HKO{1}\xrightarrow[\hspace{.0cm}]{\eta_*} \HKO{0} \xrightarrow[\hspace{.0cm}]{2_*} \HKO{0} 
\xrightarrow[\hspace{.0cm}]{\rho_*} \HKR{0} \xrightarrow[\hspace{.0cm}]{\delta_*} \HKO{1}
\xrightarrow[\hspace{.0cm}]{\rho_*} \HKU{1}
\]
\end{thm}

In the above diagram there are 20 spaces,  
$KU$ and $K(1)$ are 2-periodic and
$KO$ and $KR(1)$ are 8 periodic.
We know the homology of 12 of them.
It is the other 8 associated with $KR(1)$ that we are interested in.
Counting the suspension
maps, there are 98 maps to evaluate, 48 of them involving the $KR(1)$ spaces.
For each map, there is a spectral sequence, and 56 of them involve the $KR(1)$ spaces.
It is not necessary to know all of them to get our main results,
but it was often helpful.  
Because this information is perhaps more interesting than the main 
theorems, it has been written
up as an appendix.
Once you know the homology of all the spaces and also know the maps, it is
fairly easy to figure out how all the spectral sequences behave.
Also,  the long exact sequences
of homotopy groups have been put in the appendix as well.
In the first part of the paper we state, compute, and use, only what we need, and we
assume results not involving the $\KR{i}$.

The spaces $\KR{i}$ have been around for a long time.  When I tried to find
a reference for the homotopy groups, the experts informed me that they were
known in the 1960s to Mahowald and that there wasn't a reference because everyone
already knew them.
What might be new is that the $\KR{i}$ are also the real first Morava $K(1)$-theory.
This comes from the work of Hu-Kriz, \cite{HK}, where they compute the homotopy of
all of the real Morava $K(n)$, $KR(n)$.
This project got started because I thought $KR(2)$ would be interesting but that
I should quickly take a look at $KR(1)$ first.
From the point of view of personal satisfaction, the homology, $\HKR{6}$, was
both the most difficult to compute and the most interesting.
In the beginning, motivation was easy.  I was hoping to find something interesting.  
After the fact, it isn't clear how to motivate.

In Section \ref{known} we give the homology of the spaces that are known already
as well as state the details of Theorem \ref{exactness}.
In Section \ref{ha-ss} we state the spectral sequences we will use and
discuss how Hopf algebras help us with our computations.
After that, each section is just the computation of some $\HKR{i}$.
They are somewhat in order except except that to do $\HKR{6}$, we need to
have $\HKR{7}$ first, which is computed from $\HKR{0}$.
We finish off with an appendix computing all 98 relevant maps and spectral
sequences.

\section{Connecting to known results}
\label{known}

Our preferred generators for $\HKU{*}$ and $\HKO{*}$ come from
Hopf rings.  They are given elegant descriptions in \cite{CS}.
In   \cite[Section 25]{NituER2}, there is an alternative Hopf ring
description for $\HKO{*}$ and one can read off that for $\HKU{*}$
from \cite{RW:HR}. We will not write down these descriptions in
this paper.  It is enough to know they have nice Hopf ring names.
In the case of $\HKR{*}$, we will not get 
Hopf ring names because $KR(1)$ is not a ring spectrum.

We give the descriptions of the homologies we need in this paper.

\begin{thm} The homology of the connected component of $\KU{i}$ is as follows.
The index $k$ runs over all $k > 0$.
\[
\begin{array}{cll}
i=0 \qquad &
 \otimes_k P(x_{2k})  & \qquad V(x_{4k}) = x_{2k}  \\
i=1 \qquad & 
\otimes_k E(x_{2k+1})  &
\end{array}
\]
\end{thm}

\begin{thm} The homology of the connected component of $\KO{i}$ is as follows.
The index $k$ runs over all $k > 0$.
\[
\begin{array}{cll}
i=0 \qquad &
\otimes_k P(x_k) &  \qquad V(x_{2k}) = x_k  \\
i=1 \qquad &
 \otimes_k P(x_{2k+1}) & \\
i=2 \qquad &
\otimes_k P(x_{4k+2})& \\
i=3 \qquad &
\otimes_k E(x_{4k+3})& \\
i=4 \qquad &
\otimes_k P(x_{4k}) &  \qquad V(x_{8k}) = x_{4k}  \\
i= 5 \qquad &
\otimes_k E(x_{4k+1})&  \\
i=6 \qquad &
\otimes_k E(x_{2k}) &  \qquad V(x_{4k}) = x_{2k}  \\
i=7 \qquad &
 \otimes_k E(x_{k}) &  \qquad V(x_{2k}) = x_{k}  
\end{array}
\]
\end{thm}

\begin{thm} The homology of the connected component of $\K{i}$ is as follows.
The index $k$ runs over all $k > 0$.
\[
\begin{array}{lc}
i=0 \qquad &
TP_4(x_{4k+3}) \otimes_k E(y_{4k}) \otimes_k E(z_{8k+2}) \\
 & \\
&
 V(y_{8k}) = y_{4k}  \qquad 
V(y_{16k+4}) = z_{8k+2}
\qquad
V(y_{16k+12}) = (x_{4k+3})^2 \\
 & \\
i=1 \qquad &
E(x_{4k+1}) \otimes_k P(y_{4k+2}) \qquad V(y_{8k+2}) = x_{4k+1}
\end{array}
\]
\end{thm}

In the paper this is from, \cite{WSW:MK}, we computed $H_*(\ul{K(n)}_{\;*})$ for all $n$ and
all primes.  Slight adjustments had to be made all along the way for 
$p=2$, and it seems that they weren't all made.

In the paper, we write:
\[
\HK{0} \simeq E(x_{4k+3}) \otimes_k E(x_{2k}) 
\]
but we missed the extension $x_{4k+3}^2 = x_{8k+6}$.  So, what
is in the paper is an associated graded version.  
When the spectral sequence there is used to compute $\HK{1}$,
deep down in the gruesome depths of the paper there is a $d_1$, so
the resulting answer is correct.  
Explicitly, what it shows in that paper is that we need (in the notation
of the paper) $(e_1 a_{(0)})^2 = b_{(0)}b_{(1)}$.  
The rest follows from Hopf ring considerations
as our generators there all have nice Hopf ring names.
Something
similar happens for $K(n)$ in that paper, but again, only for $p=2$.

We can now use these results to connect to our new results.

\begin{thm}
The exactness at the middle term of $\HKO{i}\ra \HKR{i} \ra \HKO{i+1}$ of
Theorem \ref{exactness} is given explicitly as follows where, if not described,
the element maps to zero.
The index $j$ runs over all $j > 0$.
\begin{alignat*}{2}
 i&=0 \qquad &
\otimes_j P(y_j) 
\xrightarrow{y_j \ra z_j}
E(z_j)&\otimes_j P(zz_{4j+2})
\xrightarrow{zz_{4j+2} \ra (w_{2j+1})^2}
\otimes_j P(w_{2j+1}) \\
 i&=1  \qquad &
\otimes_j P(y_{2j+1}) 
\xrightarrow{y_{2j+1} \ra z_{2j+1}}
P(z_{2j+1})&\otimes_j P(zz_{4j+2})
\xrightarrow{zz_{4j+2} \ra w_{4j+2}}
\otimes_j P(w_{4j+2}) \\
i&=2 \qquad &
\otimes_j P(y_{4j+2}) 
\xrightarrow[y_{8j+6} \ra (zz_{4j+3})^2]{y_{8j+2} \ra z_{8j+2}}
P(z_{8j+2})&\otimes_j P(zz_{4j+3})
\xrightarrow{zz_{4j+3} \ra w_{4j+3}}
\otimes_j E(w_{4j+3}) \\
i&=3 \qquad &
\otimes_j E(y_{4j+3}) 
\xrightarrow{y_{8j+3} \ra z_{8j+3}}
E(z_{8j+3})&\otimes_j P(zz_{8j+4})
\xrightarrow{zz_{8j+4} \ra w_{8j+4}}
\otimes_j P(w_{4j}) \\
i&=4 \qquad &
\otimes_j P(y_{4j}) 
\xrightarrow{y_{4j} \ra z_{4j}}
E(z_{4j})&\otimes_j E(zz_{8j+5})
\xrightarrow{zz_{8j+5} \ra w_{8j+5}}
\otimes_j E(w_{4j+1}) \\
i&=5 \qquad &
\otimes_j E(y_{4j+1}) 
\xrightarrow{y_{4j+1} \ra z_{4j+1}}
E(z_{4j+1})&\otimes_j E(zz_{2j})
\xrightarrow{zz_{2j} \ra w_{2j}}
\otimes_j E(w_{2j}) \\
i&=6  \qquad &
\otimes_j E(y_{2j}) 
\xrightarrow{y_{2j} \ra (z_{j})^2}
T&P_4(z_{j})
\xrightarrow{z_{j} \ra w_{j}}
\otimes_j E(w_{j}) \\
 i&=7 \qquad &
\otimes_j E(y_{j}) 
\xrightarrow{y_{2j} \ra z_{2j}}
E(z_{2j}) &\otimes_j P(zz_{2j+1})
\xrightarrow{zz_{2j+1} \ra w_{2j+1}}
\otimes_j P(w_{j}) 
\end{alignat*}
%\[
%\begin{array}{lc}
%{\bf i=0} &
%P(y_j) 
%\xrightarrow{y_j \ra z_j}
%E(z_j)\otimes P(zz_{4j+2})
%\xrightarrow[z_j \ra 0]{zz_{4j+2} \ra (w_{2j+1})^2}
%P(w_{2j+1}) \\
%{\bf i=1} &
%P(y_{2j+1}) 
%\xrightarrow{y_{2j+1} \ra z_{2j+1}}
%P(z_{2j+1})\otimes P(zz_{4j+2})
%\xrightarrow[z_{2j+1} \ra 0]{zz_{4j+2} \ra w_{4j+2}}
%P(w_{4j+2}) \\
%{\bf i=2} \quad  &
%P(y_{4j+2}) 
%\xrightarrow[y_{8j+6} \ra (zz_{4j+3})^2]{y_{8j+2} \ra z_{8j+2}}
%P(z_{8j+2})\otimes P(zz_{4j+3})
%\xrightarrow[z_{8j+2} \ra 0,\quad (zz_{4j+3})^2 \ra 0]{zz_{4j+3} \ra w_{4j+3}}
%E(w_{4j+3}) \\
%{\bf i=3} &
%E(y_{4j+3}) 
%\xrightarrow[y_{8j+7} \ra 0]{y_{8j+3} \ra z_{8j+3}}
%E(z_{8j+3})\otimes P(zz_{4j+4})
%\xrightarrow[z_{8j+3} \ra 0 ]{zz_{8j+4} \ra w_{8j+4}}
%P(w_{4j}) \\
%{\bf i=4} &
%P(y_{4j}) 
%\xrightarrow[(y_{4j})^2 \ra 0]{y_{4j} \ra z_{4j}}
%E(z_{4j})\otimes E(zz_{8j+5})
%\xrightarrow[z_{8j+1} \ra 0 ]{zz_{8j+5} \ra w_{8j+5}}
%E(w_{4j+1}) \\
%{\bf i=5} &
%E(y_{4j+1}) 
%\xrightarrow{y_{4j+1} \ra z_{4j+1}}
%E(z_{4j+1})\otimes E(zz_{2j})
%\xrightarrow[z_{4j+1} \ra 0 ]{zz_{2j} \ra w_{2j}}
%E(w_{2j}) \\
%{\bf i=6} &
%E(y_{2j}) 
%\xrightarrow{y_{2j} \ra (z_{j})^2}
%TP_4(z_{j})
%\xrightarrow[(z_{j})^2 \ra 0 ]{z_{j} \ra w_{j}}
%E(w_{j}) \\
%{\bf i=7} &
%E(y_{j}) 
%\xrightarrow[y_{2j+1}\ra 0]{y_{2j} \ra z_{2j}}
%E(z_{2j}) \otimes P(zz_{2j+1})
%\xrightarrow[z_{2j} \ra 0 ]{zz_{2j+1} \ra w_{2j+1}}
%P(w_{j}) \\
%\end{array}
%\]
\end{thm}

\begin{remark}
The long exact sequence for $i=0$ of theorem \ref{exactness} consists
of the above maps spliced together with well-understood maps that we will
see throughout the paper.
\end{remark}

\section{Hopf algebras, fibrations, and spectral sequences}
\label{ha-ss}

We need two spectral sequences.
The homology version we use computes the homology of
a base space from the homologies of the fiber and
the total space. It is in \cite[Theorems 2.2 and 3.1]{Moore-bar}.
I think of it as the bar spectral sequence, but it should
perhaps be called the Moore spectral sequence.
Unfortunately, Moore doesn't indulge appropriately
with Hopf algebras as he clearly could have.  
Rothenberg-Steenrod,
\cite{RS}, really bring in the Hopf algebras, but neglected to
do the more general case where the total space isn't contractible.
Everyone seems to think they can do it just slightly extending
Rothenberg-Steenrod's proof, except those who think it is already in
Rothenberg-Steenrod.
The cohomology version computes the cohomology of the fiber
from the cohomologies of the base space and the total space.
This seems to originate with Eilenberg and Moore in \cite{EilMoore}.
However, my favorite reference here is \cite{Smith:LEMSS} because
this is where I learned to compute with Hopf algebras in these
spectral sequences.

We will state the two spectral sequences for the record and
then discuss the use of Hopf algebras in their computations.

\begin{prop}
	Let $F \ra E \ra B$ be a fibration of infinite loop spaces and maps.
	\begin{enumerate}
		\item
	There is a first quadrant homology spectral sequence of Hopf algebras
	\[
		E_{*,*}^2 = \Tor^{H_*(F)}_{*,*}(H_*(E),\Zq) => H_*(B).
		\]
			with $d_r \colon E_{u,v} \ra E_{u-r,v+r-1}.$
\item
	There is a second quadrant cohomology spectral sequence of Hopf algebras
	\[
		E^{*,*}_2 = \Tor_{H^*(B)}^{*,*}(H^*(E),\Zq) => H^*(F).
		\]
			with $d_r \colon E^{u,v} \ra E^{u+r,v-r+1}$
	\end{enumerate}
\end{prop}

{\bf Discussion of Hopf algebras, Tor, and differentials.}

Combining the above spectral sequences with Hopf algebras makes for a
powerful tool.  We will only discuss the homology version but everything
carries over to the cohomology version.
The general reference for Hopf algebras is \cite{MM}, but my computational
reference is \cite{Smith:LEMSS}.

We work with mod 2 homology throughout.
The Borel structure theorem (see \cite{MM}) for our graded Hopf algebras over $\Zq$ is
that they are the tensor products of algebras of the form $P(x_i)$ (polynomial),
$E(x_i)$ (exterior), and $TP_{2^j}(x_i) = P(x_i)/(x_i^{2^j})$ (truncated polynomial).
(Recall our notation is that $x_i$ is of degree $i$.)
Sub-Hopf algebras of polynomial algebras must also be polynomial.
In our Hopf algebras, we have $2_* = FV = VF$ where $F $ is the 
Frobenius (i.e. $x \ra x^2$) and $V$ is the 
	Vershiebung (i.e. the dual of the Frobenius on cohomology). 
The Hopf algebra $\Gamma[x_i]$is dual to $P(y_i)$ with $y_i$ primitive.
As such, it is $\Zq$-free on elements $\gamma_k(x_i)$ in degree $ki$.
As an algebra, it is an exterior on the generators $\gamma_{2^j}(x_i)$
of degree $2^j i$.  
We have $V(\gamma_{2^{j+1}}(x_i))  = \gamma_{2^j}(x_i)$.

There are a number of situations that arise frequently in our computations.
For example, we might find that we have an associated graded object that
is $\otimes_i E(x_i)$, but we know that when the extensions
are solved it  must be polynomial.  This becomes $\otimes_i P(x_{2i+1})$
for degree reasons.
Similarly $\otimes_i E(x_{2i})$ and $\otimes_i E(x_{4i})$, if they are really polynomial
algebras, become $\otimes_i P(x_{4i+2})$ and $\otimes_i P(x_{8i+4})$. 
If we have $\otimes_i \Gamma[x_i]$ as an associated graded object for what we know
is polynomial, we get $\otimes_i P(x_i)$.

On the other hand, if there are no extension problems, as algebras, we have
$\otimes_i \Gamma[x_{2i+1}]$
is just $\otimes_i E(y_i)$, 
and
$\otimes_i \Gamma[x_{4i+2}]$
is just $\otimes_i E(y_{2i})$.

If we have the differential algebra $E(x_i) \otimes P(y_{i+1})$ with $d^1(x_i) = y_{i+1}$,
we know that the homology in positive degrees is zero.  We will often be confronted with
the dual of this situation where we have $E(x_i) \otimes \Gamma[y_{i+1}]$ with
$d_1(y_{i+1}) = x_i$.  Again, our homology here is zero in positive degrees.
It is not always that simple though.  It often happens that we have 
$E(x_{2i+1}) \otimes \Gamma[y_{i+1}]$ and have $d_2(\gamma_2(y_{i+1})) = x_{2i+1}$.
This leaves $E(y_{i+1})$ as its homology.  When this happens, we will abuse
notation and write 
\[
E(x_{2i+1}) \otimes \Gamma[y_{i+1}] \simeq E(x_{2i+1}) \otimes E(y_{i+1})
\otimes \Gamma[y_{2i+2}]
\]
so we can see the differential and results more clearly.  This is just the associated
graded object we get from the short exact sequence of Hopf algebras
\[
E(y_{i+1}) \lra \Gamma[y_{i+1}] \lra \Gamma[y_{2i+2}]
\]
where we have written $\gamma_2(y_{i+1}) = y_{2i+2}$.
Similarly, worse happens and we need 
\[
E(x_{4i+3}) \otimes \Gamma[y_{i+1}] \simeq E(x_{4i+3}) \otimes E(y_{i+1})
\otimes E(y_{2i+2})
\otimes \Gamma[y_{4i+4}]
\]
where we have a differential taking $y_{4i+4}$ to $x_{4i+3}$ leaving
only 
$  E(y_{i+1}) 
\otimes E(y_{2i+2})
$
but with $V(y_{2i+2}) = y_{i+1}$.

To deal with our spectral sequences, we must be able to evaluate Tor.
The simple case of $\Tor_0$ is the Hopf algebra cokernel of the map
$H_*(F) \ra H_*(E)$.  There are no differentials on this zero filtration
and what remains after differentials hit it is a sub-Hopf algebra of
$H_*(B)$, i.e. the image of $H_*(E) \ra H_*(B)$.

We have a few facts to accumulate.

	\begin{enumerate}
	\item
		If $A$ is the Hopf algebra kernel of 
			the map $H_*(F) \ra H_*(E)$,
			then the higher filtrations are given by
			$\Tor^A(\Zq,\Zq)$.
		\item
			Tor commutes with tensor products. 
		\item
			$\Tor^{E(x_i)}(\Zq,\Zq) = \Gamma[y_{i+1}]$
			with $y_{i+1} $ in bidegree $(1,i)$ and
			$\gamma_{2^j}(y_{i+1})$ in bidegree $2^j$ times this.
		\item
			$\Tor^{P(x_i)}(\Zq,\Zq) = E(y_{i+1})$
			with $y_{i+1}$ in bidegree $(1,i)$.
		\item
			$\Tor^{TP_{2^k}(x_i)}(\Zq,\Zq) 
= E(y_{i+1})\otimes \Gamma[z_{2^k i + 2}]$
			with $y_{i+1}$ in bidegree $(1,i)$ and $z_{2^k i +2}$
			in bidegree $(2,2^k i)$.
		\item
			Elements in filtrations zero and one are permanent cycles.
	\end{enumerate}

If the kernel, $A$, is trivial, the spectral sequence collapses and the cokernel
is $H_*(B)$, giving us a short exact sequence $H_*(F) \ra H_*(E) \ra H_*(B)$.

Since the kernel, $A$, is a Hopf algebra, Borel's theorem applies and the above
allows us to compute Tor completely.  Differentials must start on the second or higher
filtration and they must take generators to primitives.  The primitives all live
in filtrations 0 or 1 and all generators in filtrations 2 or higher are of even 
degree.  Thus the targets of differentials must be odd degree elements in 
filtrations 0 or 1.  A fact that we will often use is that {\bf any even degree element
in filtrations 0 or 1 must survive.}

There is one more special case we need to discuss.
If we have a short exact sequence $E(y_{2i}) \ra TP_4(x_i) \ra E(x_i)$ that
takes $y_{2i}$ to $(x_i)^2$, we can compute Tor of $TP_4(x_i)$ as above and
get $E(z_{i+1})\otimes \Gamma[w_{4i+2}]$.  If we didn't know there was the
square $x_i^2 = y_{2i}$ in the middle term, but thought the middle term might
be $E(y_{2i})\otimes E(x_i)$, then Tor would be $\Gamma[z_{i+1}] \otimes
\Gamma[u_{2i+1}]$.  If we had a reason to know that this was not correct, then
a $d_1(\gamma_2(z_{i+1})) = u_{2i+1}$ would leave us with the correct answer.

\section{$\HKR{0}$}

We begin with the spectral sequence for
$$
\xymatrix{
	\KO{0} \ar[r]^2 & \KO{0} \ar[r] & \KR{0}
	}
$$
Computing $2_*$ is easy, we have 
\[
2_*(x_{2i}) = FV(x_{2i}) = F(x_i) = (x_i)^2
\]
We can read off the cokernel as $\otimes_i E(x_i)$ and the kernel as $\otimes_i P(x_{2i+1})$.
Computing Tor on the kernel we get $\otimes_j E(y_{2j})$.  Since all of these
generators are in filtrations zero and one, the spectral sequence collapses.
What we know at this stage is that we have 
\[
\otimes_i E(x_i) \subset \HKR{0} \qquad \qquad V(x_{2i}) = x_i
\]
with quotient having an associated graded object, $\otimes_j E(y_{2j})$.

We move now to a different spectral sequence, the one for
\[
	\KO{0} \lra \KR{0} \lra \KO{1}
\]
We have computed the image of $\HKO{0} \ra \HKR{0}$.  It is just
$\otimes_i E(x_i)$.  The cokernel is the object with associated 
graded object $\otimes_i E(y_{2i})$ above.  That is our zero filtration for
this spectral sequence.  
The generators are all in even degrees and so must survive.
This is all of the zeroth filtration and the zero filtration
must be a sub-Hopf algebra of $\HKO{1}$ which is polynomial
so the cokernel must be polynomial, and for degree reasons, this
must be $\otimes_j P(y_{4j+2})$.  This splits as algebras and coalgebras
and so completes our computation.

\section{$\HKR{1}$}

We start with the spectral sequence for the fibration
\[
\xymatrix{
	\KO{1} \ar[r]^2 & \KO{1} \ar[r] & \KR{1}
	}
\]
The map $2_*$ is zero because all of the generators $x_{2i+1}$ for
$\HKO{1}$ are primitive, so $V(x_{2i+1}) = 0$ giving $2_*(x_{2i+1})
= FV(x_{2i+1}) = 0$.  The cokernel is $\HKO{1} = \otimes_i P(x_{2i+1})$ and
so is the kernel.  We now know the zeroth filtration and taking
Tor of the kernel, we get exterior generators $y_{2i}$ in filtration 1.  
The
spectral sequence collapses because all the generators are in
filtrations 0 and 1.  We still have extension problems though.
Again, we move to the next spectral sequence for
\[
	\KO{1} \lra \KR{1} \lra \KO{2}
\]
We have computed the image of $\HKO{1} \ra \HKR{1}$.  It is just
$\otimes_i P(x_{2i+1})$.  There is no kernel, 
so the spectral sequence collapses and is just the 
cokernel in the zeroth filtration.
This becomes a short
exact sequence of Hopf algebras.
\[
		\HKO{1} \lra \HKR{1} \lra \HKO{2} \
\]
But this is just
\[
 \otimes_i P(x_{2i+1}) \lra \HKR{1} \lra \otimes_i P(y_{4i+2}) 
\]
and so splits as algebras, giving us most of our answer.  
There is an extension problem to solve to get $V(y_{4i+2}) = x_{2i+1}$.

For that we use the spectral sequence for:
\[
	\KR{0} \lra \text{ * }  \lra  \KR{1}
\]
Computing Tor of $\HKR{0}$ we get 
\[
\Gamma[w_{k}] \otimes_k E(ww_{4k+3})
\]
We should note that we have to use the $\Zq$ in degree zero for $H_0(\KR{0})$
to get the $w_{1}$ above.

This is way too big.  Remember, we know the answer here.  To get this down to
size, we must take the first possible differential, i.e. we must have
$d_3(\gamma_4(w_{k})) \not = 0$.  This element is degree $4k$ (in
the fourth filtration) so the
differential hits an element in the first filtration in degree $4k-1$.
There are two possibilities, but it must hit one of them, and we don't
need to know which just yet.
All that is left after these differentials is an exterior algebra,
$\otimes_i E(z_i)$ with generators in filtration 1 and an exterior algebra
$\otimes_i E(\gamma_2(w_i))$ with generators in filtration 2.  This is precisely
the size of the known answer so these differentials must indeed happen.

We know that the answer is polynomial, so the Frobenius must be
injective.  The Frobenius cannot raise filtration so the injective
Frobenius on the first filtration gives us $\otimes_i P(z_{2i+1})$, forcing
(to get the correct answer) the Frobenius to inject on the second
filtration to get $\otimes_i P(w_{4i+2})$.  The only ambiguity in the first
filtration is about which elements in degrees $4i+3$ have survived.
We know that the element in degree $2i+1$ in the
first filtration must square to the element in degree $4i+2$,
and this is unambiguously $x_{4i+2} = V( \gamma_2(x_{4i+2}))$.
But we know that we must have $\gamma_2(x_{4i+2}) =
(\gamma_2(x_{2i+1}))^2 $ because of the injectivity of $F$.
But now we have just computed $VF$ on $\gamma_2(x_{2i+1})$ and found
it non-zero.  Consequently, $VF=FV$ must also be non-zero so
that $V$ is non-zero.  We get our result that $V$ of every generator
of $P(y_{4i+2})$ is a generator of $P(x_{2i+1})$ as desired.

\section{$\HKR{2}$}

We start with the spectral sequence for
\[
\xymatrix{
	\KO{2} \ar[r]^2 & \KO{2} \ar[r] & \KR{2}
	}
\]
The map $2_*$ is zero because all of the generators $x_{4i+2}$ for
$\HKO{2}$ are primitive, so $V(x_{4i+2}) = 0$ giving $2_*(x_{4i+2})
= FV(x_{4i+2}) = 0$.  The cokernel is $\HKO{2} = \otimes_i P(x_{4i+2})$ and
so is the kernel.  We now know the zeroth filtration and taking
Tor of the kernel, we get $\otimes_i E(y_{4i+3})$.  The
spectral sequence collapses because all the generators are in
filtrations 0 and 1.  We still have extension problems though.

What we have from the spectral sequence is the short exact sequence:
\[
\otimes_i  P(x_{4i+2}) \lra \HKR{2} \lra \otimes_i E(y_{4i+3}) 
\]
There is an extension problem
we need to solve, namely $(y_{4i+3})^2$ from filtration 1 is $x_{8i+6}$
in filtration 0.  Once this is done, we would have the algebra structure.

To solve this problem we look at the spectral sequence for
\[
\KR{2} \llra{\eta} \KR{1} \lra \K{1}
\]
We have maps
\[
\HKR{2} \lra P(z_{2i+1})\otimes_i P(zz_{4i+2}) \lra E(w_{4i+1})\otimes_i P(ww_{4i+2})
\]
Our calculation so far shows that $\HKR{2}$ is generated by primitives.
We we know from our computation of $\HKR{1}$
that $V(zz_{4i+2}) = z_{2i+1}$, so all the primitives in $\HKR{1}$ are
in $\otimes_i P(z_{2i+1})$.  
Since primitives map to primitives, we see that $\otimes_i P(zz_{4i+2})$ is in
the cokernel.  It is even degree and in filtration zero so is a sub-algebra
of $\HK{1}$, so it must be our $\otimes_i P(ww_{4i+2})$ in our known answer.
This accounts for all of the even degree generators and squares in
$\HK{1}$.

If the element $y_{4i+3}$ from $\HKR{2}$ is in the kernel, then Tor will
give rise to an element in filtration 1 of degree $4i+4$.  This element
would have to survive, but we have all of the even degree generators
and squares we need, so $y_{4i+3}$ maps to $z_{4i+3}$ because it is
the only primitive in that degree.  However, $z_{4i+3}$ is a polynomial
generator, so $y_{4i+3}$ must also be a polynomial generator, solving
our extension problem.

We can go one step further.  If $x_{8i+2}$ doesn't map to $(z_{4i+1})^2$,
this last element would be even degree in the cokernel where we don't
need any more even degree elements, so it does map accordingly.
We get a rare short exact
sequence.
\[
 \HKR{2} \lra \HKR{1} \lra \HK{1} .
\]

\section{$\HKR{3}$}

We start with the spectral sequence for 
\[
	\KR{2} \lra \text{ * }  \lra \KR{3}
\]
Since $\HKR{2} \simeq P(x_{8i+2})\otimes_i P(x_{4i+3})$, computing Tor is
easy, it is just
\[
E(x_{8i+3}) \otimes_i E(x_{4i})
\]
and since the generators are all in filtration 1, it collapses.
All we have left are extension problems.

Next we use the spectral sequence for:
\[
\xymatrix{
	\KO{3} \ar[r]^2 & \KO{3} \ar[r] & \KR{3}
	}
\]
The homology, $\HKO{3}$, is generated by primitives, so $2_*$ is zero.
We get the cokernel is $\otimes_i E(x_{4i+3})$ and so is the kernel. 
The $E^2$ term of the spectral sequence is
\[
E(x_{4i+3}) \otimes_i \Gamma[y_{4i}]
\]

This  is much too big compared with our first spectral sequence.
The only way to cut it down to the right size is with
\[
d_2(\gamma_2(y_{4i})) = x_{8i-1}
\]
This leaves $E(x_{8i+3})\otimes_i E(y_{4i})$ as with the first one,
but now we know that the $E(x_{8i+3})$ is the image of $\HKO{3}$
in $\HKR{3}$ and the cokernel has an associated graded object
of $\otimes_i E(y_{4i})$.

We can move on to the spectral sequence for
\[
	\KO{3} \lra  \KR{3} \lra \KO{4}
\]
We just computed the cokernel.  It is even degree in filtration zero
and all of the elements must survive.  Since this cokernel is a
subalgebra of the polynomial algebra $\HKO{4}$, this solves
all of our extension problems giving $(y_{4i})^2 = y_{8i}$ so
we have the expected polynomial algebra $\otimes_i P(y_{8i+4})$, completing
our computation.

\section{$\HKR{4}$}

We use the spectral sequence coming from
\[
\xymatrix{
	\KO{4} \ar[r]^2 & \KO{4} \ar[r] & \KR{4}
	}
\]
As in the $\KR{0}$ case, $\HKO{4}$ is bipolynomial.
The cokernel is just $\otimes_i E(x_{4i})$ and the kernel is
$P(x_{8i+4})$. We take Tor of this to get exterior generators $y_{8i+5}$
in filtration 1.  Since all our generators are in filtrations zero and one,
the spectral sequence collapses.  For purely degree reasons, there can
be no extension problems given that we know $\otimes_i E(x_{4i})$ is a sub-algebra.

\section{$\HKR{5}$}

We use the spectral sequence for the fibration
\[
	\KR{4} \lra \text{ * }  \lra \KR{5}
\]
Since $\HKR{4} \simeq E(x_{4k})\otimes_k E(y_{8k+5})$, Tor is
\[
\Gamma[z_{4k+1}] \otimes_k \Gamma[w_{8k+6}]
\]
The only possible sources for differentials are in degrees divisible by 4,
but the only odd degree elements are in degree 1 mod 4, so there can be
no differentials.  Furthermore, there are no algebra extension problems.
In filtration one there are only elements $z_{4k+1}$ and $w_{8k+6}$,
so there is nothing for them to square to.  In filtration two, the elements
are in degrees $8k+2$ and $16k+12$ and again there are no elements in 
filtrations 0 or 1 to square to.  Continue inductively on filtration.
The degrees never work out to have extensions.
This spectral sequence gives a complete description of $V$ as well.

\section{$\HKR{7}$}

Note that we have skipped $\HKR{6}$.  It is the hardest to compute
and all our previous techniques failed us.  We need $\HKR{7}$
to solve the problems with $\HKR{6}$.

We use the cohomology spectral sequence for the fibration 
\[
\xymatrix{
\KR{7} \ar[r] & \KO{0} \ar[r]^2 & \KO{0}
}
\]
The homology of $\KO{0}$ is bipolynomial with $\HKO{0} \simeq \otimes_i P(x_i)$
and $V(x_{2i}) = x_i$.  So we get $H^*(\KO{0})$ is the same.  Evaluating
$2^*$, takes $2^*(x_{2i}) = FV(x_{2i}) = F(x_i) = (x_i)^2$.  The cokernel
is $\otimes_i E(x_i)$ with $V$ as before.  Since $V(x_{2i+1}) = 0$, the kernel
is $\otimes_i P(x_{2i+1})$.  Tor of the kernel is $\otimes_i E(w_{2i})$ with generators
in the first filtration.  Since all of the generators are in the first
2 filtrations, the spectral sequence collapses.  Since we know the $V$
on filtration zero ($\otimes_i E(x_i)$), we can dualize and we get the homology has
$\otimes_i P(y_{2i+1})$ in it and there is 
the $\otimes_i E(ww_{2i})$ (dual to $\otimes_i E(w_{2i})$) 
as well, but it could
have extension problems we need to solve.

To show that the $\otimes_i E(ww_{2i})$ really is exterior, we take a quick
look at the homology spectral sequence for
\[
\KO{7} \llra{2} \KO{7} \lra \KR{7}
\]
The first map is zero because $\HKO{7} \simeq \otimes_i E(z_i)$ and $F$ is zero, so the
cokernel contains $\otimes_i E(z_{2i})$ in the zero filtration and this
subalgebra must survive, we now have the desired exterior subalgebra.

\section{$\HKR{6}$ }

This is both the hardest to compute and the most interesting.
We start with our usual fibration
\[
\xymatrix{
	\KO{6} \ar[r]^2 & \KO{6} \ar[r] & \KR{6}
	}
\]
$\HKO{6} \simeq \otimes_k E(x_{2k})$ so $2_* = VF = 0 $ because $F$ is zero.
That means all of $\HKO{6}$ is the cokernel and it all survives
because it is even degree.
Using our second spectral sequence for
\[
\KO{6} \lra \KR{6} \lra \KO{7}
\]
we know that the first map injects so there is no kernel.
The spectral sequence collapses with $\HKO{7}$ the cokernel of the map.
This gives the short exact sequence
\[
\xymatrix{
	 \otimes_i E(x_{2i})	 \ar[r] & \HKR{6} \ar[r] & \otimes_i E(y_i) 
	}
\]
The goal here is to solve the extension problem $(y_i)^2 = x_{2i}$.
We do already know
that $V(x_{4i}) = x_{2i}$ on the first part and $V(y_{2i}) = y_i$
on the second part.

We use the spectral sequence
\[
\KR{6} \lra * \lra \KR{7}
\]
to prove our result.  Observe that Tor of $E(x_{2i})\otimes_i E(y_i)$ is
\[
\Gamma[w_{2i+1}] \otimes_i \Gamma[ww_i]
= 
E(w_{2i+1})\otimes_i \Gamma[w_{4i+2}] \otimes_i E(ww_i) \otimes_i \Gamma[ww_{2i}]
\]
and that Tor of $\otimes_i TP_4(y_i)$ is
\[
 \Gamma[w_{4i+2}] \otimes_i E(ww_i) 
\]
If the extension exists, there is a $d_1(ww_{2i}) = w_{2i-1}$.
We don't need $d_1$ for our result though.  Note that no matter what the
extension is, in Tor we have
\[
 \Gamma[w_{4i+2}] \otimes_i E(ww_i) 
\]
We rewrite  this just a bit as
\[
 \Gamma[w_{4i+2}] \otimes_i E(ww_{2i}) \otimes_i E(ww_{2i+1}).
\]
Note that this is precisely the correct size for our known result of $\HKR{7}$.
That doesn't prove our result yet though.  We do know that any even
degree element in filtration 1 or 2 must survive, and so we know we
have $E(w_{4i+2}) \otimes_i E(ww_{2i})$ already no matter what extensions
there are.  
One of $w_{4i+2}$ or $ww_{4i+2}$ must be exterior and the other must square.
The only thing to square to is $ww_{8i+4}$.  Because we know the answer
and all these elements must survive, this must be part of the polynomial
part of the answer, so we must have $(ww_{4i})^2 = ww_{8i}$.
It doesn't really matter which of the elements is exterior.
What we know from this is that we have all the elements we need in
degrees $4i+2$ that are generators, primitives, or squares.

If we had a case where $(y_{2i})^2 = 0$ in $\HKR{6}$, then from the above
discussion, we would have Tor giving us a $\Gamma[z_{2i+1}] =
E(z_{2i+1}) \otimes 
\Gamma[z_{4i+2}]$ 
We would not have the 
$\Gamma[z_{4i+2}]$ 
unless this happens.
The $z_{4i+2}$ is in filtration 2 so must survive, but we already
have enough elements in this degree, so this cannot happen.

We now know that $(y_{2i})^2 = x_{4i}$ always.  We have
\[
x_{2i} = V(x_{4i}) = VF(y_{2i}) = FV(y_{2i}) = F(y_i) = (y_i)^2
\]
This solves the extension problem for all $y_i$ with $i$ odd.

\begin{center}
{\LARGE \bf Appendix}
\end{center}

\section{Homologies}

To make this appendix somewhat more self contained, we write down
all of the homologies we need.
In the first part of the paper we write $\otimes_i P(x_i)$ to be clear that we
are taking the polynomial algebra on generators $x_i$ for all $i$.
The tensor products clutter up the notation.  That kind of precision
isn't necessary in this appendix, so here, when we write $P(x_i)$,
we mean $\otimes_i P(x_i)$.  The tensor product is understood.
This simplifies the notation significantly.

\begin{thm} The homology of the connected component of $\KR{i}$ is as follows.
	If the Vershiebung isn't described, it is zero.
\begin{alignat*}{4}
&i=0\hspace{.5in} 
&
E(x_k) &\otimes P(y_{4k+2})
&
\hspace{.25in}V(x_{2k}) &= x_k
& \\
&i=1\hspace{.5in}
&
P(x_{2k+1})&\otimes P(y_{4k+2}) 
&
\hspace{.55in}V(y_{4k+2}) &= x_{2k+1}
&\\
&i=2\hspace{.5in}
&
P(x_{8k+2}) &\otimes P(y_{4k+3})
&
&
& \\
&i=3
&
E(x_{8k+3}) &\otimes P(y_{8k+4})
&
&
& \\
&i=4
&
E(x_{4k}) &\otimes E(y_{8k+5})
&
\hspace{.25in}V(x_{8k}) &= x_{4k}
& \\
&i=5
&
E(x_{4k+1}) &\otimes E(y_{2k})
&
\hspace{.25in}V(y_{4k}) &= y_{2k}
\qquad V(y_{8k+2}) &=  x_{4k+1}
\\
&i=6
&
T&P_4(x_{k})
&
\hspace{.25in}V(x_{2k}) &= x_{k}
& \\
&i=7
&
E(x_{2k})
&\otimes 
P(y_{2k+1}) 
&
\hspace{.25in}V(x_{4k}) &= x_{2k}
&
\end{alignat*}
\end{thm}

\begin{thm} The homology of the connected component of $\KU{i}$ is as follows.
\[
\begin{array}{cll}
i=0 \qquad &
 P(x_{2i})  & \qquad V(x_{4i}) = x_{2i}  \\
i=1 \qquad & 
E(x_{2i+1})  &
\end{array}
\]
\end{thm}

\begin{thm} The homology of the connected component of $\KO{i}$ is as follows.
\[
\begin{array}{cll}
i=0 \qquad &
P(x_i) &  \qquad V(x_{2i}) = x_i  \\
i=1 \qquad &
 P(x_{2i+1}) & \\
i=2 \qquad &
P(x_{4i+2})& \\
i=3 \qquad &
E(x_{4i+3})& \\
i=4 \qquad &
P(x_{4i}) &  \qquad V(x_{8i}) = x_{4i}  \\
i= 5 \qquad &
E(x_{4i+1})&  \\
i=6 \qquad &
E(x_{2i}) &  \qquad V(x_{4i}) = x_{2i}  \\
i=7 \qquad &
 E(x_{i}) &  \qquad V(x_{2i}) = x_{i}  
\end{array}
\]
\end{thm}

\begin{thm} The homology of the connected component of $\K{i}$ is as follows.
\[
\begin{array}{lc}
i=0 \qquad &
TP_4(x_{4k+3}) \otimes E(y_{4k}) \otimes E(z_{8k+2}) \\
 & \\
&
 V(y_{8k}) = y_{4k}  \qquad 
V(y_{16k+4}) = z_{8k+2}
\qquad
V(y_{16k+12}) = (x_{4k+3})^2 \\
 & \\
i=1 \qquad &
E(x_{4k+1}) \otimes P(y_{4k+2}) \qquad V(y_{8k+2}) = x_{4k+1}
\end{array}
\]
\end{thm}

\section{Introduction and ground rules}

We will not reproduce computations done in the first part of 
the paper,
hereafter
referred to as the "main paper."  Previous to the main paper, the maps, homologies,
and spectral sequences of the spaces $\KU{i}$, $\KO{i}$, and $\K{i}$ were all
known.  We will not recompute these but only describe them.  Many of the maps
and spectral sequences associated with $\KR{i}$ are not computed in the main
paper.  However, relying on the main paper, we do know $\HKR{i}$ as well as
the homologies for all of the previously known spaces.  
Some of the maps and spectral sequences have already been computed in the
main paper, but not, by any means, all.  When we have a new map or spectral
sequence we will do more than just describe it, but we will give details of
the proof of what is new.
For every spectral 
sequence we study here, we know the answer, which is often quite helpful.
In fact, often the argument for a differential or the solution to an
extension problem is "because we know the answer."  Rather than keep
repeating this phrase, we will just assert differentials and extension
problem solutions if they must come about "because we know the answer."

Because there are 98 spectral sequences, we developed self-explanatory notation
for them so we can refer to them if necessary.
The general form is a sequence of fibrations
\[
\ul{X}_{\;i}
\lra
\ul{Y}_{\;i}
\lra
\ul{Z}_{\;i}
\lra
\ul{X}_{\;i+1}
\lra
\ul{Y}_{\;i+1}
\lra
\ul{Z}_{\;i+1}
\lra
\ul{X}_{\;i+2}
\lra \cdots
\]
First we compute the map 
$H_*( \ul{X}_{\;i}) 
\lra
H_*( \ul{Y}_{\;i}) 
$
and from that we compute the spectral sequence for 
$H_*( \ul{Z}_{\;i}) $.
With the exception of one spectral sequence out of 98, this gives the
map
$H_*( \ul{Y}_{\;i}) 
\lra
H_*( \ul{Z}_{\;i}) .
$
When we move to the next spectral sequence, i.e. for 
$\ul{Y}_{\;i}
\lra
\ul{Z}_{\;i}
\lra
\ul{X}_{\;i+1}$,
we already know the first map and we can repeat the story and move
on to the next.
Knowing the first map of each spectral sequence and the answer makes
most of them quite easy to deal with.

Of the 4 such sequences we describe, only one requires some work
to start with the first map, namely $\HKR{1} \llra{\eta} \HKR{0}$.
The only second map that doesn't come out of the spectral sequence is
the map for $\HK{5} \ra \HKR{7}$ from the spectral sequence {\bf RKR557}.
The problem is solved in the very next spectral sequence, { \bf KRR576}.

We have also labeled the infrequent short exact sequences so the curious
can find them easily.

\section{$\HKO{i} => \HKO{i+1} \quad $ {\bf OOi(i+1)}}

We use the bar spectral sequence for
\[
\KO{i} \lra * \lra \KO{i+1}
\]

{\bf i=0, OO01}
\[
\HKO{0} = P(x_{i}) => \HKO{1} = P(y_{2i+1})
\]
Tor of $P(x_i)$ is $E(y_i)$ with $y_i$ in filtration 1.  Solving all extension
problems, $(y_i)^2 = y_{2i}$, 
 gives $P(y_{2i+1})$.

{\bf i=1, OO12}
\[
\HKO{1} = P(x_{2i+1}) => \HKO{2} = P(y_{4i+2})
\]
Tor of $P(x_{2i+1})$ is $E(y_{2i})$ with $y_{2i}$ in filtration 1.  Solving all extension
problems, 
$(y_{2i})^2 = y_{4i}$,
gives $P(y_{4i+2})$.

{\bf i=2, OO23}
\[
\HKO{2} = P(x_{4i+2}) => \HKO{3} = E(y_{4i+3})
\]
Tor of $P(x_{4i+2})$ is $E(y_{4i+3})$ with $y_{4i+3}$ in filtration 1.  
This gives $E(y_{4i+3})$.

{\bf i=3, OO34}
\[
\HKO{3} = E(x_{4i+3}) => \HKO{4} = P(y_{4i})
\]
Tor of $E(x_{4i+3})$ is $\Gamma[y_{4i}]$ with $y_{4i}$ in filtration 1.  
We have $(y_{4i})^2 = y_{8i}$ and corresponding formulas on the $\gamma$s.
This gives $P(y_{4i})$.

{\bf i=4, OO45}
\[
\HKO{4} = P(x_{4i}) => \HKO{5} = E(y_{4i+1})
\]
Tor of $P(x_{4i})$ is $E(y_{4i+1})$ with $y_{4i+1}$ in filtration 1.  
This gives $E(y_{4i+1})$.

{\bf i=5, OO56}
\[
\HKO{5} = E(x_{4i+1}) => \HKO{6} = E(y_{2i})
\]
Tor of $E(x_{4i+1})$ is $\Gamma[y_{4i+2}]= E(y_{2i})$.

{\bf i=6, OO67}
\[
\HKO{6} = E(x_{2i}) => \HKO{7} = E(y_{i})
\]
Tor of $E(x_{2i})$ is $\Gamma[y_{2i+1}]= E(y_i)$.

{\bf i=7, OO70}
\[
\HKO{7} = E(x_{i}) => \HKO{0} = P(y_{i})
\]
Tor of $E(x_{i})$ is $\Gamma[y_{i}]$ with $y_{i}$ in filtration 1.  
We have $(y_{i})^2 = y_{2i}$ and corresponding formulas on the $\gamma$ generators.
This gives $P(y_{i})$.

\section{$\HKU{i} => \HKU{i+1} \quad $ {\bf UUi(i+1)}}

We use the bar spectral sequence for
\[
\KU{i} \lra * \lra \KU{i+1}
\]

{\bf i=0, UU01}
\[
\HKU{0} = P(x_{2i}) => \HKU{1} = E(y_{2i+1})
\]
Tor of $P(x_{2i})$ is $E(y_{2i+1})$ with $y_{2i+1}$ in filtration 1.  
This  gives $E(y_{2i+1})$.

{\bf i=1, UU10}
\[
\HKU{1} = E(x_{2i+1}) => \HKU{0} = P(y_{2i})
\]
Tor of $E(x_{2i+1})$ is $\Gamma[y_{2i}]$ with $y_{2i}$ in filtration 1.  Solving all extension
problems, $(y_{2i})^2 = y_{4i}$ and similar formulas on the $\gamma$ generators,  gives $P(y_{2i})$.

\section{$\HK{i} => \HK{i+1} \quad $ {\bf KKi(i+1)}}

We use the bar spectral sequence for
\[
\K{i} \lra * \lra \K{i+1}
\]

{\bf i=0, KK01}
\[
\HK{0} = TP_4(x_{4i+3})\otimes E(y_{4i})\otimes   E(z_{8i+2}) => 
\HK{1} = E(x_{4i+1}) \otimes P(y_{4i+2})
\]
Tor of 
\[
TP_4(x_{4i+3})\otimes E(y_{4i})\otimes   E(z_{8i+2}) 
\] 
is 
\[
E(x_{4i}) \otimes \Gamma[x_{16i+14}] \otimes 
\Gamma[y_{4i+1}]
\otimes
\Gamma[z_{8i+3}]
\] 
with $x_{4i}$, $y_{4i+1}$, and $z_{8i+3}$ in filtration 1  
and $x_{16i+14}$ in filtration 2.
We rewrite 
$\Gamma[y_{4i+1}]$
as
$E(y_{4i+1}) \otimes 
\Gamma[y_{8i+2}]$
with the $y_{8i+2}$ in filtration 2.
We rewrite 
$\Gamma[z_{8i+3}]$
as
$E(z_{8i+3}) \otimes 
\Gamma[z_{16i+6}]$
with the $y_{16i+6}$ in filtration 2.
Our Tor is now 
\[
E(x_{4i})\otimes E(x_{4i+1}) \otimes E(x_{8i+3}) \otimes \Gamma[z_{4i+2}]
\]
where the exterior generators are in filtration 1 and the $\Gamma$ generator
is in filtration 2.
Now we rewrite $\Gamma[z_{4i+2}]$ as $E(z_{4i+2})\otimes \Gamma[z_{8i+4}]$
with the $\Gamma$ generator in filtration 4.  We now have a $d_3$,
\[
d_3(z_{8i+4}) = x_{8i+3}.
\]
We are left with 
\[
E(x_{4i})\otimes E(x_{4i+1}) \otimes  E(z_{4i+2})
\]
Solve the extension problems to get the known result
$\HK{1} = E(x_{4i+1}) \otimes P(y_{4i+2})$.
Those solutions are $(x_{4i})^2 = x_{8i}$ and $(z_{4i+2})^2 = x_{8i+4}$.

{\bf i=1, KK10}
\[
\HK{1} = E(x_{4i+1}) \otimes P(y_{4i+2}) =>
\HK{0} = TP_4(x_{4i+3})\otimes E(y_{4i})\otimes   E(z_{8i+2}) 
\]
Tor is 
\[
\Gamma[x_{4i+2}]\otimes E(y_{4i+3})
\]
Rewrite $\Gamma[x_{4i+2}] $ as 
\[
\Gamma[x_{8i+2}] 
\otimes
E(x_{8i+6}) \otimes \Gamma[x_{16i+12}]
\]
There is one extension problem, $(y_{4i+3})^2 = x_{8i+6}$. 
After this we have \[
TP_4(y_{4i+3}) \otimes \Gamma[x_{8i+2}] \otimes \Gamma[x_{16i+12}]
\]
Continuing to rewrite, this is
\[
TP_4(y_{4i+3}) \otimes E(x_{8i+2}) \otimes \Gamma[x_{8i+4}]
\]
and the $\Gamma[x_{8i+4}]$ is $E(x_{4i})$.

\section{The sequence $\quad \KO{1} \llra{\eta} \KO{0} \ra \KU{0} 
\ra \KO{2} \llra{\eta} \KO{1} \ra \cdots$}

\begin{itemize}
\setlength\itemsep{1em}
\item{$\KO{1} \llra{\eta} \KO{0} \ra \KU{0}\qquad $ {\bf OOU100}}
\[
P(x_{2i+1}) \xrightarrow{x_{2i+1} \ra y_{2i+1} } P(y_i) 
\xrightarrow[y_{2i+1} \ra 0]{y_{2i} \ra z_{2i}}
P(z_{2i})
\]
The spectral sequence is just a short exact sequence
\[
\text{\bf S.E.S.} \qquad  \HKO{1} \llra{\eta_*} \HKO{0} \lra \HKU{0}.
\]

\item{$ \KO{0} \ra \KU{0} \ra \KO{2}\qquad $ {\bf OUO002}}
\[
P(x_i) \xrightarrow[x_{2i+1} \ra 0]{x_{2i} \ra y_{2i} } P(y_{2i})
\xrightarrow{zero} P(z_{4i+2})
\]
There is nothing in the cokernel.  The kernel is $P(x_{2i+1})$ so
Tor is $E(w_{2i})$.  Solving the extension problems, $(w_{2i})^2 = w_{4i}$, 
gives $P(z_{4i+2})$.

\item{$  \KU{0} \ra \KO{2}\llra{\eta} \KO{1} \qquad $ {\bf UOO021}}
\[
P(x_{2i})
\xrightarrow{zero} P(y_{4i+2})
\xrightarrow{y_{4i+2} \ra (z_{2i+1})^2}
P(z_{2i+1})
\]
The cokernel is $P(y_{4i+2})$ in filtration 0 and the kernel is $P(x_{2i})$.  Tor
on this is $E(w_{2i+1})$ with generators in filtration 1.  
We get $(w_{2i+1})^2 = y_{4i+2}$. 

\item{$  \KO{2}\llra{\eta} \KO{1} \ra \KU{1}\qquad $ {\bf OOU211}}
\[
P(x_{4i+2})
\xrightarrow{x_{4i+2} \ra (y_{2i+1})^2}
P(y_{2i+1})
\xrightarrow[(y_{2i+1})^2 \ra 0 ]{y_{2i+1} \ra z_{2i+1}}
E(z_{2i+1})
\]
We get a short exact sequence
\[
\text{\bf S.E.S.} \qquad  \HKO{2} \lra \HKO{1} \lra \HKU{1}.
\]

\item{$   \KO{1} \ra \KU{1}\ra \KO{3}\qquad $ {\bf OUO113}}
\[
P(x_{2i+1})
\xrightarrow[(x_{2i+1})^2 \ra 0 ]{x_{2i+1} \ra y_{2i+1}}
E(y_{2i+1})
\xrightarrow{zero} 
E(z_{4i+3})
\]
The cokernel is zero and the kernel is $P((x_{2i+1})^2)$.  Tor of
this is $E(w_{4i+3})$, our answer.

\item{\bf $    \KU{1}\ra \KO{3} \llra{\eta} \KO{2}\qquad$ {\bf UOO132}}
\[
E(x_{2i+1})
\xrightarrow{zero} 
E(y_{4i+3})
\xrightarrow{zero} 
P(z_{4i+2})
\]
The cokernel in filtration zero is $E(y_{4i+3})$ and the kernel is $E(x_{2i+1})$.
Tor of this is $\Gamma[w_{2i}]$.  There is a differential
\[
d_2(\gamma(w_{2i})) = y_{4i-1}
\]
leaving only $E(w_{2i})$.  Solving the extension problems, $(w_{2i})^2 = w_{4i}$, 
gives $P(z_{4i+2})$.

\item{\bf $     \KO{3} \llra{\eta} \KO{2} \ra \KU{2}\qquad $ {\bf OOU322}}
\[
E(x_{4i+3})
\xrightarrow{zero} 
P(y_{4i+2})
\xrightarrow{y_{4i+2} \ra z_{4i+2}} 
P(z_{2i})
\]
The cokernel in filtration zero is $P(z_{4i+2})$.  The kernel is $E(x_{4i+3})$.
Tor of this is $\Gamma[w_{4i}]$ starting in filtration 1.
There are no differentials and everything in $\Gamma[w_{4i}]$ squares 
non-trivially to something in the same filtration, i.e. $(w_{4i})^2 = w_{8i}$
and all of the corresponding $\gamma$'s do the same, giving $P(z_{4i})$.  

\item{ $      \KO{2} \ra \KU{2}\ra \KO{4} \qquad $ {\bf OUO224}}
\[
P(x_{4i+2})
\xrightarrow{x_{4i+2} \ra y_{4i+2}} 
P(y_{2i})
\xrightarrow[y_{4i+2} \ra 0]{y_{4i} \ra z_{4i}}
P(z_{4i})
\]
We get a short exact sequence
\[
\text{\bf S.E.S.} \qquad  \HKO{2} \lra \HKU{2} \lra \HKO{4}.
\]

\item{ $ \KU{2}\ra \KO{4} \llra{\eta} \KO{3} \qquad $ {\bf UOO243}}
\[
P(x_{2i})
\xrightarrow[x_{4i+2} \ra 0]{x_{4i} \ra y_{4i}}
P(y_{4i})
\xrightarrow{zero} 
E(z_{4i+3})
\]
There is no cokernel.  The kernel is $P(x_{4i+2})$.  Tor on this
is $E(w_{4i+3})$ and we are done.

\item{ $  \KO{4} \llra{\eta} \KO{3} \ra \KU{3} \qquad $ {\bf OOU433}}
\[
P(x_{4i})
\xrightarrow{zero} 
E(y_{4i+3})
\xrightarrow{y_{4i+3} \ra z_{4i+3}}
E(z_{2i+1})
\]
The cokernel is 
$E(y_{4i+3})$
and the kernel is 
$P(x_{4i})$.
Tor on this is
$E(y_{4i+1})$ and we have our answer.

\item{ $   \KO{3} \ra \KU{3} \ra \KO{5} \qquad $ {\bf OUO335}}
\[
E(x_{4i+3})
\xrightarrow{x_{4i+3} \ra y_{4i+3}}
E(y_{2i+1})
\xrightarrow[y_{4i+3} \ra 0]{y_{4i+1} \ra z_{4i+1}}
E(z_{4i+1})
\]
This gives a short exact sequence
\[
\text{\bf S.E.S.} \qquad  \HKO{3} \lra \HKU{3} \lra \HKO{5}.
\] 

\item{ $  \KU{3} \ra \KO{5} \llra{\eta} \KO{4} \qquad $ {\bf UOO354}}
\[
E(x_{2i+1})
\xrightarrow[x_{4i+3} \ra 0]{x_{4i+1} \ra y_{4i+1}}
E(y_{4i+1})
\xrightarrow{zero} 
P(z_{4i})
\]
There is no cokernel and the kernel is $E(x_{4i+3})$.  Tor
of this is $\Gamma[w_{4i}]$.  Solving extensions, $(w_{4i})^2 = w_{8i}$ and
corresponding formulas on the $\gamma$ generators,  gives
$P(z_{4i})$.

\item{ $  \KO{5} \llra{\eta} \KO{4}\ra \KU{4} \qquad $ { \bf OOU544}}
\[
E(x_{4i+1})
\xrightarrow{zero} 
P(y_{4i})
\xrightarrow{y_{4i} \ra (z_{2i})^2 }
P(z_{2i})
\]
The cokernel is $P(y_{4i})$ and the kernel 
$E(x_{4i+1})$.
Tor of this is 
$\Gamma[w_{4i+2}]=E(w_{2i})$.
We have $(w_{2i})^2 = y_{4i}$.

\item{ $  \KO{4}\ra  \KU{4}\ra  \KO{6}  \qquad $ {\bf OUO446}}
\[
P(x_{4i})
\xrightarrow{x_{4i} \ra (y_{2i})^2 }
P(y_{2i})
\xrightarrow[(y_{2i})^2 \ra 0 ]{y_{2i} \ra z_{2i}} 
E(z_{2i})
\]
We get a short exact sequence
\[
\text{\bf S.E.S.} \qquad  \HKO{4} \lra \HKU{4} \lra \HKO{6}.
\]

\item{ $   \KU{4}\ra  \KO{6} \llra{\eta} \KO{5}  \qquad $ {\bf UOO465}}
\[
P(x_{2i})
\xrightarrow[(x_{2i})^2 \ra 0 ]{x_{2i} \ra y_{2i}} 
E(y_{2i})
\xrightarrow{zero} 
E(z_{4i+1})
\]
There is no cokernel.  The kernel is $P((x_{2i})^2)$.  The
Tor of this is $E(w_{4i+1})$, our answer.

\item{ $     \KO{6} \llra{\eta} \KO{5} \ra \KU{5}  \qquad $ {\bf OOU655}}
\[
E(x_{2i})
\xrightarrow{zero} 
E(y_{4i+1})
\xrightarrow{zero} 
E(z_{2i+1})
\]
The cokernel is $E(y_{4i+1})$ and the kernel is $E(x_{2i})$.  Tor
of this is $\Gamma[w_{2i+1}]$.  We must have a differential
\[
d_2(\gamma_2(w_{2i+1})) = y_{4i+1}
\]
All that is left is the $E(w_{2i+1})$, our answer.

\item{ $  \KO{5} \ra \KU{5} \ra \KO{7}  \qquad $ {\bf OUO557}}
\[
E(x_{4i+1})
\xrightarrow{zero} 
E(y_{2i+1})
\xrightarrow{y_{2i+1} \ra z_{2i+1}} 
E(z_i)
\]
The cokernel is 
$E(y_{2i+1})$
and the kernel is
$E(x_{4i+1})$.
Tor of this is $\Gamma[w_{4i+2}]= E(w_{2i})$.  

\item{ $   \KU{5} \ra \KO{7} \llra{\eta} \KO{6} \qquad $ {\bf UOU576}}
\[
E(x_{2i+1})
\xrightarrow{x_{2i+1} \ra y_{2i+1}} 
E(y_i)
\xrightarrow[y_{2i+1} \ra 0]{y_{2i} \ra z_{2i}} 
E(z_{2i})
\]
We get a short exact sequence
\[
\text{\bf S.E.S.} \qquad  \HKU{5} \lra \HKO{7} \lra \HKO{6}.
\]

\item{ $    \KO{7} \llra{\eta} \KO{6} \ra \KU{6} \qquad $ {\bf OOU766}}
\[
E(x_i)
\xrightarrow[x_{2i+1} \ra 0]{x_{2i} \ra y_{2i}} 
E(y_{2i})
\xrightarrow{zero} 
P(z_{2i})
\]
There is no cokernel.  The kernel is $E(x_{2i+1})$.  Tor of this is $\Gamma[w_{2i}]$.
Solving the extension problems we get our answer, $P(z_{2i})$.

\item{ $  \KO{6}\ra \KU{6} \ra \KO{0}  \qquad $ {\bf OUO660}}
\[
E(x_{2i})
\xrightarrow{zero} 
P(y_{2i})
\xrightarrow{y_{2i} \ra (z_{i})^2} 
P(z_{i})
\]
The cokernel is $P(y_{2i})$ and the kernel is $E(x_{2i})$.  Tor
of this is $\Gamma[w_{2i+1}] = E(w_{i})$.  The 
extension problem is solved by $(w_{i})^2 = y_{2i}$.

\item{ $  \KU{6} \ra \KO{0}\llra{\eta} \KO{7}  \qquad $ {\bf UOO607}}
\[
P(x_{2i})
\xrightarrow{x_{2i} \ra (y_{i})^2} 
P(y_{i})
\xrightarrow[(y_i)^2 \ra 0 ]{y_{i} \ra z_{i}} 
E(z_i)
\]
We get a short exact
sequence
\[
\text{\bf S.E.S.} \qquad  \HKU{6} \lra \HKO{0} \ra \HKO{7}.
\]

\item{ $  \KO{0}\llra{\eta} \KO{7} \ra \KU{7}  \qquad $ {\bf OOU077}}
\[
P(x_{i})
\xrightarrow[(x_i)^2 \ra 0 ]{x_{i} \ra y_{i}} 
E(y_i)
\xrightarrow{zero} 
E(z_{2i+1})
\]
The cokernel is zero and the kernel is $P((x_i)^2)$.  Tor
of this is $E(w_{2i+1})$.

\item{ $   \KO{7} \ra \KU{7} \ra \KO{1}  \qquad $ {\bf OUO771}}
\[
E(x_i)
\xrightarrow{zero} 
E(y_{2i+1})
\xrightarrow{zero} 
P(z_{2i+1})
\]
The cokernel is $E(y_{2i+1})$ and the kernel is $E(x_i)$.  Tor
of this is $\Gamma[w_{i}]$. 
There is a differential, $d_2(\gamma_2(w_i)) = y_{2i-1}$ leaving
$E(w_i)$.  Solving extensions, $(w_i)^2 = w_{2i}$, we get our $P(w_{2i+1})$.

\item{ $    \KU{7} \ra \KO{1} \llra{\eta} \KO{0}  \qquad $ {\bf UOO710}}
\[
E(x_{2i+1})
\xrightarrow{zero} 
P(y_{2i+1})
\xrightarrow{y_{2i+1} \ra z_{2i+1}} 
P(z_i)
\]
The cokernel is $P(z_{2i+1})$ and the kernel is $E(x_{2i+1})$.  Tor
of this is $\Gamma[w_{2i}]$. 
To get our answer we must have $(w_{2i})^2 = w_{4i}$ in filtration
1 and $FV = VF$ will give us the squares in all the higher filtrations
ending with $P(z_{2i})$.

\end{itemize}

\section{The sequence $\quad \KU{0} \llra{2} \KU{0} \ra \K{0} \ra \KU{1} 
\llra{2} \KU{1} \ra \cdots$}

\begin{itemize}
\setlength\itemsep{1em}

\item{$\KU{0} \llra{2} \KU{0} \ra \K{0} \qquad $ { \bf UUK000}}
\[
P(x_{2i}) 
\xrightarrow[x_{4i+2} \ra 0 ]{x_{4i} \ra (y_{2i})^2} 
P(y_{2i})
\xrightarrow[y_{8i+6} \ra (zz_{4i+3})^2]{y_{4i} \ra z_{4i},\quad {y_{8i+2}\ra z_{8i+2}} }
TP_4(zz_{4i+3})\otimes E(z_{4i})\otimes E(z_{8i+2})
\]
The cokernel is $E(y_{2i})$ and the kernel is $P(x_{4i+2})$.  Tor of
this is $E(w_{4i+3})$.  We have only the one extension problem,
$(w_{4i+3})^2 = x_{8i+6}$.

\item{$ \KU{0} \ra \K{0}\ra \KU{1} \qquad $ { \bf UKU001}}
\[
P(x_{2i})
\xrightarrow[x_{8i+6} \ra (yy_{4i+3})^2]{x_{4i} \ra y_{4i},\quad {x_{8i+2}\ra y_{8i+2}} }
\]
\[
TP_4(yy_{4i+3})\otimes E(y_{4i})\otimes E(y_{8i+2})
\xrightarrow[
y_{8i+2} \ra 0, \quad (yy_{4i+3})^2 \ra 0
]{
y_{4i} \ra 0, \quad
yy_{4i+3} \ra z_{4i+3}}
E(z_{2i+1})
\]
The cokernel is $E(yy_{4i+3})$ and the kernel is $P((x_{2i})^2)$.
Tor of this is $E(w_{4i+1})$.

\item{$  \K{0}\ra \KU{1} \llra{2} \KU{1} \qquad $ {\bf KUU011}}
\[
TP_4(xx_{4i+3})\otimes E(x_{4i})\otimes E(x_{8i+2})
\]
\[
\xrightarrow[
x_{8i+2} \ra 0, \quad (xx_{4i+3})^2 \ra 0
]{
x_{4i} \ra 0, \quad
xx_{4i+3} \ra y_{4i+3}}
E(y_{2i+1})
\xrightarrow{zero}
E(z_{2i+1})
\]
The cokernel is $E(y_{4i+1})$ and the kernel $E(x_{2i})$.
Tor of the kernel
is
$\Gamma[w_{2i+1}] \simeq E(w_{2i+1})\otimes \Gamma[w_{4i+2}]$.
We get a differential
\[
d_2(ww_{4i+2}) = y_{4i+1}
\]
All that is left is $E(w_{2i+1})$.

\item{$   \KU{1} \llra{2} \KU{1} \ra \K{1}  \qquad $ {\bf UUK111}}
\[
E(x_{2i+1})
\xrightarrow{zero}
E(y_{2i+1})
\xrightarrow[y_{4i+3} \ra 0 ]{y_{4i+1} \ra z_{4i+1}}
E(z_{4i+1}) \otimes P(zz_{4i+2}).
\]
The cokernel is $E(y_{2i+1})$ and the kernel is $E(x_{2i+1})$.
Tor of this is $\Gamma[w_{2i}]$.
We have a differential
\[
d_2(\gamma_2(w_{2i})) = y_{4i-1}
\]
We are left with $E(y_{4i+1})$ in filtration zero, 
and $E(w_{2i})$ with generators in filtration 1.  This last 
all have squares, $(w_{2i})^2 = w_{4i}$, giving $P(w_{4i+2})$.

\item{$   \KU{1} \ra \K{1} \ra \KU{0}  \qquad $ {\bf UKU110}}
\[
E(x_{2i+1})
\xrightarrow[x_{4i+3} \ra 0 ]{x_{4i+1} \ra y_{4i+1}}
E(y_{4i+1}) \otimes P(yy_{4i+2})
\xrightarrow[y_{4i+1} \ra 0 ]{yy_{4i+2} \ra z_{4i+2}}
P(z_{2i})
\]
The cokernel is $P(yy_{4i+2})$ and the kernel is $E(x_{4i+3})$.
Tor of this is $\Gamma[w_{4i}]$.  Squaring everything in 
$\Gamma$ gives $P(w_{4i})$.

\item{$   \K{1} \ra \KU{0} \llra{2} \KU{0}  \qquad $ {\bf KUU100}}
\[
E(x_{4i+1}) \otimes P(xx_{4i+2})
\xrightarrow[x_{4i+1} \ra 0 ]{xx_{4i+2} \ra y_{4i+2}}
P(y_{2i})
\xrightarrow[y_{4i+2} \ra 0 ]{y_{4i} \ra (z_{2i})^2} 
P(z_{2i})
\]
The cokernel is $P(y_{4i})$ and the kernel is $E(x_{4i+1})$.
Tor of this is $\Gamma[w_{4i+2}] = E(w_{2i})$.
We have $(w_{2i})^2 = y_{4i}$ giving our answer.

\end{itemize}

\section{The sequence $\quad \KO{0} \llra{2} \KO{0} \ra \KR{0} \ra \KO{1} 
\llra{2} \KO{1} \ra \cdots$}

\begin{itemize}
\setlength\itemsep{1em}

\item{$\KO{0} \llra{2} \KO{0} \ra \KR{0} \qquad $ {\bf OOR000  }}
\[
P(x_{i}) 
\xrightarrow[x_{2i+1} \ra 0 ]{x_{2i} \ra (y_{i})^2} 
P(y_{i})
\xrightarrow[(y_{i})^2 \ra 0]{y_{i} \ra z_{i} }
E(z_i)\otimes P(zz_{4i+2})
\]
The cokernel is $E(y_i)$ and the kernel is $P(x_{2i+1})$.  Tor
of this is $E(w_{2i})$. 
We have $(w_{2i})^2 = w_{4i}$ giving 
 $P(w_{4i+2})$.

\item{$\KO{0} \ra \KR{0}\ra \KO{1} \qquad $ {\bf ORO001 }}
\[
P(x_{i})
\xrightarrow[(x_{i})^2 \ra 0]{x_{i} \ra y_{i} }
E(y_i)\otimes P(yy_{4i+2})
\xrightarrow[y_i \ra 0]{yy_{4i+2} \ra (z_{2i+1})^2 }
P(z_{2i+1})
\]
The cokernel is $P(yy_{4i+2})$ 
The kernel is
$P((x_i)^2)$.  Tor of this is $E(w_{2i+1})$.  
We have  $(w_{2i+1})^2 = yy_{4i+2}$.

\item{$ \KR{0}\ra \KO{1}\llra{2} \KO{1} \qquad  $ {\bf ROO011 }}
\[
E(x_i)\otimes P(xx_{4i+2})
\xrightarrow[x_i \ra 0 ]{xx_{4i+2} \ra (y_{2i+1})^2 }
P(y_{2i+1})
\xrightarrow{zero}
P(z_{2i+1})
\]
The cokernel is $E(y_{2i+1})$ and the kernel is $E(x_i)$.
Tor of this is $\Gamma[w_i]$.
We have differentials
\[
d_2(\gamma_2(w_i) ) = y_{2i-1}.
\]
We are left with $E(w_i)$.  
We have $(w_i)^2 = w_{2i}$ so we get
 $P(w_{2i+1})$.

\item{$  \KO{1}\llra{2} \KO{1} \ra \KR{1} \qquad $ {\bf OOR111 }}
\[
P(x_{2i+1})
\xrightarrow{zero}
P(y_{2i+1})
\xrightarrow{y_{2i+1} \ra z_{2i+1}}
P(z_{2i+1}) \otimes P(zz_{4i+2})
\]
The cokernel is $P(y_{2i+1})$ and the kernel is $P(x_{2i+1})$.
Tor of this is $E(w_{2i})$ and
We have $(w_{2i})^2 = w_{4i}$ so we get
$P(w_{4i+2})$.

\item{$   \KO{1} \ra \KR{1}\ra \KO{2} \qquad $ {\bf ORO112 }}
\[
P(x_{2i+1})
\xrightarrow{x_{2i+1} \ra y_{2i+1}}
P(y_{2i+1}) \otimes P(yy_{4i+2})
\xrightarrow[y_{2i+1} \ra 0]{yy_{4i+2} \ra z_{4i+2} }
P(z_{4i+2})
\]
There is no kernel so this is a short exact sequence
\[
\text{\bf S.E.S.} \qquad  \HKO{1} \lra \HKR{1} \lra \HKO{2}.
\]

\item{$    \KR{1} \ra \KO{2} \llra{2} \KO{2} \qquad $ { \bf ROO122 }}
\[
P(x_{2i+1}) \otimes P(xx_{4i+2})
\xrightarrow[x_{2i+1} \ra 0]{xx_{4i+2} \ra y_{4i+2} }
P(y_{4i+2})
\xrightarrow{zero}
P(z_{4i+2})
\]
The cokernel is zero and the kernel is $P(x_{2i+1})$.  Tor
of this is $E(w_{2i})$.  
We have $(w_{2i})^2 = w_{4i}$ so we get
$P(w_{4i+2})$.

\item{$ \KO{2} \llra{2} \KO{2} \ra \KR{2} \qquad $ { \bf OOR222 }}
\[
P(x_{4i+2})
\xrightarrow{zero}
P(y_{4i+2})
\xrightarrow[y_{8i+6} \ra (zz_{4i+3})^2]{y_{8i+2} \ra z_{8i+2}}
P(z_{8k+2})\otimes P(zz_{4i+3})
\]
The cokernel is 
$P(y_{4i+2})$
and the kernel is
$P(x_{4i+2})$.
Tor of the kernel is $E(w_{4i+3})$ and we have 
 $(w_{4i+3})^2 = x_{8i+6}$.

\item{$  \KO{2} \ra \KR{2} \ra \KO{3} \qquad $ { \bf ORO223  }}
\[
P(x_{4i+2})
\xrightarrow[x_{8i+6} \ra (yy_{4i+3})^2]{x_{8i+2} \ra y_{8i+2}}
P(y_{8k+2})\otimes P(yy_{4i+3})
\xrightarrow[y_{8i+2} \ra 0,\quad (yy_{4i+3})^2 \ra 0]{yy_{4i+3} \ra z_{4i+3}}
E(z_{4i+3})
\]
The cokernel is $E(yy_{4i+3})$ and there is no kernel.
We get a short exact sequence. 
\[
\text{\bf S.E.S.} \qquad  \HKO{2} \lra \HKR{2} \lra \HKO{3}.
\]

\item{$  \KR{2} \ra \KO{3} \llra{2} \KO{3} \qquad $ { \bf ROO233  }}
\[
P(x_{8k+2})\otimes P(xx_{4i+3})
\xrightarrow[x_{8i+2} \ra 0,\quad (xx_{4i+3})^2 \ra 0]{xx_{4i+3} \ra y_{4i+3}}
E(y_{4i+3})
\xrightarrow{zero}
E(z_{4i+3})
\]
There is no cokernel.  The kernel is $P(x_{8i+2}) \otimes P((xx_{4i+3})^2)$.
Tor of this is $E(w_{8i+3}) \otimes E(ww_{8i+7})$.

\item{$ \KO{3} \llra{2} \KO{3} \ra \KR{3} \qquad $ { \bf OOR333   }}
\[
E(x_{4i+3})
\xrightarrow{zero}
E(y_{4i+3})
\xrightarrow[y_{8i+7} \ra 0]{y_{8i+3} \ra z_{8i+3}}
E(z_{8i+3})\otimes P(zz_{8i+4})
\]
The cokernel is $E(y_{4i+3})$ and
the kernel is $E(x_{4i+3})$.  Tor of this is
$\Gamma[w_{4i}]$.
We must have a differential
\[
d_2(\gamma_2(w_{4i}) ) = y_{8i-1}
\]
This leaves $E(y_{8i+3}) \otimes E(w_{4i})$.  
We have $(w_{4i})^2 = w_{8i}$
giving $P(zz_{8i+4})$.

\item{$  \KO{3} \ra \KR{3} \ra \KO{4} \qquad $ { \bf ORO334    }}
\[
E(x_{4i+3})
\xrightarrow[x_{8i+7} \ra 0]{x_{8i+3} \ra y_{8i+3}}
E(y_{8i+3})\otimes P(yy_{8i+4})
\xrightarrow[y_{8i+3} \ra 0]{yy_{8i+4} \ra z_{8i+4}}
P(z_{4i})
\]
The cokernel
is $P(yy_{8i+4})$ and the kernel is
$E(x_{8i+7})$.  Tor of this is $\Gamma[w_{8i}]$.
We have $(w_{8i})^2 = w_{16i}$ and with $FV=VF$, this $\Gamma$ becomes
$P(z_{8i})$.

\item{$  \KR{3} \ra \KO{4} \llra{2} \KO{4} \qquad $ { \bf ROO344   }}
\[
E(x_{8i+3})\otimes P(xx_{8i+4})
\xrightarrow[x_{8i+3} \ra 0]{xx_{8i+4} \ra y_{8i+4}}
P(y_{4i})
\xrightarrow[y_{8i+4} \ra 0]{y_{8i} \ra (z_{4i})^2}
P(z_{4i})
\]
The cokernel is $P(y_{8i})$ and the kernel is $E(x_{8i+3})$.
Tor of this is $\Gamma[w_{8i+4}] = E(w_{4i})$.  
We have $(w_{4i})^2 = y_{8i}$ to get
$P(z_{4i})$.

\item{$  \KO{4} \llra{2} \KO{4} \ra \KR{4} \qquad $ { \bf OOR444   }}
\[
P(x_{4i})
\xrightarrow[x_{8i+4} \ra 0]{x_{8i} \ra (y_{4i})^2}
P(y_{4i})
\xrightarrow[(y_{4i})^2 \ra 0]{ y_{4i} \ra z_{4i}}
E(z_{4i})\otimes E(zz_{8i+5})
\]
The cokernel is $E(y_{4i})$.
The kernel is $P(x_{8i+4})$.  Tor of this is
$E(w_{8i+5})$.

\item{$ \KO{4} \ra \KR{4} \ra \KO{5}  \qquad $ { \bf ORO445    }}
\[
P(x_{4i})
\xrightarrow[(x_{4i})^2 \ra 0]{ x_{4i} \ra y_{4i}}
E(y_{4i})\otimes E(yy_{8i+5})
\xrightarrow[y_{4i} \ra 0]{ yy_{8i+5} \ra z_{8i+5}}
E(z_{4i+1})
\]
The cokernel is $E(yy_{8i+5})$.  The kernel is $P((x_{4i})^2)$.
Tor of this is $E(w_{8i+1})$.

\item{$  \KR{4} \ra \KO{5} \llra{2} \KO{5}  \qquad $ { \bf ROO455    }}
\[
E(x_{4i})\otimes E(xx_{8i+5})
\xrightarrow[x_{4i} \ra 0]{ xx_{8i+5} \ra y_{8i+5}}
E(y_{4i+1})
\xrightarrow{zero}
E(z_{4i+1})
\]
The cokernel is $E(y_{8i+1})$ and
the kernel is $E(x_{4i})$.  Tor of this is
$\Gamma[w_{4i+1}]$  
with 
\[
d_2(\gamma_2(w_{4i+1})) = y_{8i+1}
\]
What is left is $E(w_{4i+1})$.

\item{$  \KO{5} \llra{2} \KO{5} \ra \KR{5}  \qquad $ { \bf OOR555     }}
\[
E(x_{4i+1})
\xrightarrow{zero}
E(y_{4i+1})
\xrightarrow{ y_{4i+1} \ra z_{4i+1}}
E(z_{4i+1}) \otimes E(zz_{2i})
\]
The cokernel is $E(y_{4i+1})$ and the kernel is $E(x_{4i+1})$.
Tor of this is $\Gamma[w_{4i+2}] = E(w_{2i})$.

\item{$  \KO{5} \ra \KR{5} \ra \KO{6}  \qquad $ { \bf ORO556      }}
\[
E(x_{4i+1})
\xrightarrow{ x_{4i+1} \ra y_{4i+1}}
E(y_{4i+1}) \otimes E(yy_{2i})
\xrightarrow[y_{4i+1} \ra 0]{ yy_{2i} \ra z_{2i}}
E(z_{2i})
\]
The cokernel is $E(yy_{2i})$ and there is no kernel.
We get a short exact sequence 
\[
\text{\bf S.E.S.} \qquad  \HKO{5} \lra \HKR{5} \lra \HKO{6}.
\]

\item{$  \KR{5} \ra \KO{6} \llra{2} \KO{6}  \qquad $ { \bf ROO566      }}
\[
E(x_{4i+1}) \otimes E(xx_{2i})
\xrightarrow[x_{4i+1} \ra 0]{ xx_{2i} \ra y_{2i}}
E(y_{2i})
\xrightarrow{zero}
E(z_{2i})
\]
There is no cokernel.  The kernel is $E(x_{4i+1})$.
Tor of this is $\Gamma[w_{4i+2}]= E(w_{2i})$.

\item{$  \KO{6} \llra{2} \KO{6} \ra \KR{6}  \qquad $ { \bf OOR666      }}
\[
E(x_{2i})
\xrightarrow{zero}
E(y_{2i})
\xrightarrow{y_{2i} \ra (z_i)^2}
TP_4(z_i)
\]
The cokernel is $E(y_{2i})$ 
 and the kernel is $E(x_{2i})$.
Tor of this is $\Gamma[w_{2i+1}]= E(w_i)$.
We have $(w_i)^2 = y_{2i}$.

\item{$   \KO{6} \ra \KR{6} \ra \KO{7}  \qquad $ { \bf ORO667      }}
\[
E(x_{2i})
\xrightarrow{x_{2i} \ra (y_i)^2}
TP_4(y_i)
\xrightarrow{y_{i} \ra z_i}
E(z_i)
\]
The cokernel is $E(y_i)$ and there is no kernel.
We get a short exact sequence
\[
\text{\bf S.E.S.} \qquad  \HKO{6} \lra \HKR{6} \lra \HKO{7}.
\]

\item{$    \KR{6} \ra \KO{7} \llra{2} \KO{7}  \qquad $ { \bf ROO677      }}
\[
TP_4(x_i)
\xrightarrow{x_{i} \ra y_i}
E(y_i)
\xrightarrow{zero}
E(z_i)
\]
The cokernel is zero and the kernel is $E((x_i)^2)$.  Tor
of this is $\Gamma[w_{2i+1}] = E(w_i)$.

\item{$ \KO{7} \llra{2} \KO{7} \ra \KR{7}  \qquad $ { \bf OOR777       }}
\[
E(x_i)
\xrightarrow{zero}
E(y_i)
\xrightarrow[y_{2i+1} \ra 0]{y_{2i} \ra z_{2i}}
E(z_{2i})\otimes P(zz_{2i+1})
\]
The cokernel is $E(y_i)$ and the kernel is $E(x_i)$.
Tor is $\Gamma[w_i]$.  
We need
\[
d_2(\gamma_2(w_i)) = y_{2i-1}.
\]
We are left with 
$E(y_{2i})\otimes E(w_i)$. 
We must have $(w_i)^2 = w_{2i}$ to get our $P(zz_{2i+1})$.

\item{$ \KO{7} \ra \KR{7} \ra \KO{0}  \qquad $ { \bf ORO770      }}
\[
E(x_i)
\xrightarrow[x_{2i+1} \ra 0]{x_{2i} \ra y_{2i}}
E(y_{2i})\otimes P(yy_{2i+1})
\xrightarrow[y_{2i} \ra 0 ]{yy_{2i+1} \ra z_{2i+1}}
P(z_i)
\]
The cokernel is $P(yy_{2i+1})$ and the kernel is $E(x_{2i+1})$.
Tor of this is $\Gamma[w_{2i}]$.  
We have $(w_{2i})^2 = w_{4i}$, which, along with $FV=VF$ gives
$P(z_{2i})$.

\item{$ \KR{7} \ra \KO{0} \llra{2} \KO{0}  \qquad $ { \bf ROO700      }}
\[
E(x_{2i})\otimes P(xx_{2i+1})
\xrightarrow[x_{2i} \ra 0 ]{xx_{2i+1} \ra y_{2i+1}}
P(y_i)
\xrightarrow[y_{2i+1} \ra 0]{y_{2i} \ra (z_i)^2}
P(z_i)
\]
The cokernel is $P(y_{2i})$ and the kernel is $E(x_{2i})$.
Tor of this is $\Gamma[w_{2i+1}]= E(w_i)$.  
We have $(w_i)^2 =
 y_{2i}$ to get 
$P(w_i)$.

\end{itemize}

\section{$\HKR{i} => \HKR{i+1}\qquad $ {\bf RRi(i+1)}}

We use the bar spectral sequence for
\[
\KR{i} \lra * \lra \KR{i+1}
\]

{\bf i=0, RR01 }
\[
\HKR{0} = E(x_i) \otimes P(xx_{4i+2})
\lra * \lra 
\HKR{1} = P(y_{2i+1}) \otimes P(yy_{4i+2})
\]
Tor is $\Gamma[w_i] \otimes E(ww_{4i+3})$.
We have a differential
\[
d_3(\gamma_4(w_i)) = ww_{4i-1}
\]
This leaves $E(w_i)$ with generators in filtration 1 and  $ E(\gamma_2(w_i))$
with generators in filtration 2.  
We have $(w_i)^2 = w_{2i}$ giving $P(y_{2i+1})$ and
$ (\gamma_2(w_i))^2 =
 \gamma_2(w_{2i})$
giving $P(yy_{4i+2})$.

{\bf i=1, RR12 }
\[
\HKR{1} = P(x_{2i+1}) \otimes P(xx_{4i+2})
\lra * \lra 
\HKR{2} = P(y_{8i+2}) \otimes P(yy_{4i+3})
\]
Tor is $E(w_{2i})\otimes E(ww_{4i+3})$.  
We have 
$(w_{2i})^2 = w_{4i}$ and $(ww_{4i+3})^2 = w_{8i+6}$.

{\bf i=2, RR23 }
\[
\HKR{2} = P(x_{8i+2}) \otimes P(xx_{4i+3})
\lra * \lra 
\HKR{3} = E(y_{8i+3})\otimes P(yy_{8i+4})
\]
Tor is 
$ E(w_{8i+3}) \otimes E(ww_{4i})$.
We have $(ww_{4i})^2 = ww_{8i}$ giving $P(yy_{8i+4})$.

 {\bf i=3, RR34  }
\[
\HKR{3} = E(x_{8i+3})\otimes P(xx_{8i+4})
\lra * \lra 
\HKR{4} = E(y_{4i}) \otimes E(yy_{8i+5})
\]
Tor is $\Gamma[w_{8i+4}] \otimes E(ww_{8i+5})$
with 
$\Gamma[w_{8i+4}] = E(w_{4i})$.

 {\bf i=4, RR45  }
\[
\HKR{4} = E(x_{4i}) \otimes E(xx_{8i+5})
\lra * \lra 
\HKR{5} = E(y_{4i+1})\otimes E(yy_{2i})
\]
Tor is $\Gamma[w_{4i+1}]\otimes \Gamma[ww_{8i+6}]$.
Rewrite this as 
$E(w_{4i+1}) \otimes \Gamma[ww_{4i+2}]$ and
then again as our answer.

{\bf i=5, RR56   }
\[
\HKR{5} = E(x_{4i+1})\otimes E(yy_{2i})
\lra * \lra 
\HKR{6} = TP_4(y_i)
\]
Tor is $\Gamma[w_{4i+2}] \otimes \Gamma[ww_{2i+1}]$.
Rewritten, this is $E(w_{2i})$ and $E(ww_i)$.
We have $(ww_i)^2 = w_{2i}$.

{\bf i=6, RR67   }
\[
\HKR{6} = TP_4(x_i)
\lra * \lra 
\HKR{7} = E(y_{2i}) \otimes P(yy_{2i+1})
\]
Tor is $E(w_i) \otimes \Gamma[ww_{4i+2}]$.
We have $(w_i)^2 = w_{2i}$ giving $P(yy_{2i+1})$
and $\Gamma[ww_{4i+2}]$ is just $E(y_{2i})$.

{\bf i=7, RR70   }
\[
\HKR{7} = E(x_{2i}) \otimes P(xx_{2i+1})
\lra * \lra 
\HKR{0} = E(y_i)\otimes P(yy_{4i+2})
\]
Tor is $\Gamma[w_{2i+1}] \otimes E( ww_{2i})$.
We have 
$\Gamma[w_{2i+1}] $ is $E(y_i)$ and
after $(ww_{2i})^2 = w_{4i}$, we get 
$P(yy_{4i+2})$.

\section{The sequence $\quad \KR{1} \llra{\eta} \KR{0} \ra \K{0} \ra \KR{2} 
\llra{\eta} \KR{1} \ra \cdots$}

\begin{itemize}
\setlength\itemsep{1em}

\item{$\KR{1} \llra{\eta} \KR{0} \ra \K{0} \qquad $ {\bf RRK100 }}
\[
P(x_{2i+1}) \otimes P(xx_{4i+2})
\xrightarrow[(x_{2i+1})^2 \ra 0]{x_{2i+1} \ra y_{2i+1},\quad xx_{4i+2} \ra yy_{4i+2}}
E(y_{i}) \otimes P(yy_{4i+2}) 
\]\[
\xrightarrow[
y_{2i+1} \ra 0,\quad
yy_{4i+2} \ra 0]{
y_{8i+6} \ra (z_{4i+3})^2,\quad 
y_{8i+2}\ra zzz_{8i+2},\quad y_{4i}\ra zz_{4i}}
TP_4(z_{4i+3})\otimes E(zz_{4i}) \otimes E(zzz_{8i+2})
\]
We are in new territory now because we don't know the first map.  
To compute it, we use
\[
\xymatrix{
\KO{1} \ar[r] \ar[d] & \KO{0}  \ar[d] \\
\KR{1} \ar[r] \ar[d] & \KR{0}  \ar[d] \\
\KO{2} \ar[r] \ & \KO{1}   \\
}
\]
We know all of the maps in homology except the horizontal one in the middle.
We know the top horizontal map from {\bf OOU100 }
and the  bottom horizontal map from {\bf OOU211 }.
The left vertical maps are from {\bf ORO112} and the right from {\bf ORO001}.
Algebraically, we have 
\[
\xymatrix{
P(x_{2i+1}) \ar[d]^{x_{2i+1} \ra x_{2i+1}} \ar[r]^{x_{2i+1} \ra y_{2i+1}} 
& P(y_i)  \ar[d]^{y_i \ra y_i} \\
\; \; P(x_{2i+1})\otimes P(xx_{4i+2}) 
\ar[r] \ar[d]^{xx_{4i+2} \ra yy_{4i+2}} & 
\qquad  E(y_i)\otimes P(yy_{4i+2})  \ar[d]^{yy_{4i+2} \ra (yy_{2i+1})^2} \\
P(yy_{4i+2}) \ar[r]^{yy_{4i+2} \ra (yy_{2i+1})^2} \ & P(yy_{2i+1})   \\
}
\]
A diagram chase gives the first map as listed above.  

The cokernel is $E(y_{2i})$ and the kernel is $P((x_{2i+1})^2)$.
Tor of this is $E(ww_{4i+3})$.  We have
$(ww_{4i+3})^2 = y_{8i+6}$.

\item{$ \KR{0} \ra \K{0} \ra \KR{2} \qquad $ {\bf RKR002 }}
\[
E(x_{i}) \otimes P(xx_{4i+2}) 
\xrightarrow[
x_{2i+1} \ra 0,\quad
xx_{4i+2} \ra 0]{
x_{8i+6} \ra (y_{4i+3})^2,\quad 
x_{8i+2}\ra yyy_{8i+2},\quad x_{4i}\ra yy_{4i}}
\]\[
TP_4(y_{4i+3})\otimes E(yy_{4i}) \otimes E(yyy_{8i+2})
\xrightarrow{zero}
P(z_{8i+2}) \otimes P(zz_{4i+3})
\]
The cokernel is $E(y_{4i+3})$ and
the kernel is $E(x_{2i+1}) \otimes P(xx_{4i+2})$.  Tor
of this is
$\Gamma[w_{2i}] \otimes E(ww_{4i+3})$.
We need
\[
d_2( \gamma_2(w_{2i}) ) = y_{4i-1}
\]
This leaves us with 
$E(w_{2i}) \otimes E(ww_{4i+3})$.
We have $(w_{2i})^2 = w_{4i}$ 
and $(ww_{4i+3})^2 = w_{8i+6}$ to get our answer.

\item{$  \K{0} \ra \KR{2} \llra{\eta} \KR{1} \qquad $ {\bf KRR021 }}
\[
TP_4(x_{4i+3})\otimes E(xx_{4i}) \otimes E(xxx_{8i+2})
\xrightarrow{zero}
\]
\[
P(y_{8i+2}) \otimes P(yy_{4i+3})
\xrightarrow[y_{8i+2} \ra (z_{4i+1})^2]{yy_{4i+3} \ra z_{4i+3}}
P(z_{2i+1})\otimes P(zz_{4i+2}) 
\]
The cokernel 
is 
$P(y_{8i+2}) \otimes P(yy_{4i+3})$
and we know the kernel.  Tor of
the kernel, when you combine all 3 of the terms, is
$\Gamma[wwww_{4i+2}]$ starting in filtration 2, with exterior generators
in filtration 1 given by $w_{4i}$, $ww_{4i+1}$, and $www_{8i+3}$. 
We know from the main paper that the $yy_{4i+3}$ inject to $\HKR{1}$
so can't be hit by a differential.  So, we have
\[
d_3(\gamma_2(wwww_{4i+2})) = www_{8i+3}
\]
This leaves us with 
an exterior algebra with generators in filtration 1,
 $E(w_{4i}) \otimes E(ww_{4i+1})$
and 
an exterior algebra with generators in filtration 2,
$E(wwww_{4i+2})$.  
We have $(w_{4i})^2 = w_{8i}$,
$(wwww_{4i+2})^2 = w_{8i+4}$, and  $(ww_{4i+1})^2 =  y_{8i+2}$.

\item{$  \KR{2} \llra{\eta} \KR{1} \ra \K{1} \qquad $ {\bf RRK211  }}
\[
P(x_{8i+2}) \otimes P(xx_{4i+3})
\xrightarrow[x_{8i+2} \ra (y_{4i+1})^2]{xx_{4i+3} \ra y_{4i+3}}
P(y_{2i+1})\otimes P(yy_{4i+2}) 
\]\[
\xrightarrow[
y_{4i+3} \ra 0, \quad
(y_{4i+1})^2 \ra 0
]
{
yy_{4i+2} \ra zz_{4i+2}, \quad
y_{4i+1} \ra z_{4i+1}
}
E(z_{4i+1})\otimes P(zz_{4i+2}) 
\]
The cokernel is $E(y_{4i+1})\otimes P(yy_{4i+2})$ and there is no kernel.
We get a rare short exact
sequence.
\[
  \text{\bf S.E.S.} \qquad  \HKR{2} \lra \HKR{1} \lra \HK{1} .
\]

\item{$  \KR{1} \ra \K{1}\ra  \KR{3} \qquad $ {\bf RKR113   }}
\[
P(x_{2i+1})\otimes P(xx_{4i+2}) 
\xrightarrow[ 
x_{4i+3} \ra 0, \quad
(x_{4i+1})^2 \ra 0
]
{
xx_{4i+2} \ra yy_{4i+2}, \quad
x_{4i+1} \ra y_{4i+1}
}
\]
\[
E(y_{4i+1})\otimes P(yy_{4i+2}) 
\xrightarrow{zero}
E(z_{8i+3})\otimes P(zz_{8i+4})
\]
There is no cokernel and 
the kernel is $P(x_{4i+3}) \otimes P((x_{4i+1})^2)$.
Tor is $E(w_{4i}) \otimes E(ww_{8i+3})$.  
We have  $(w_{4i})^2 = w_{8i}$ giving the $P(zz_{8i+4})$.

\item{$  \K{1}\ra  \KR{3} \llra{\eta} \KR{2} \qquad $ {\bf KRR132    }}
\[
E(x_{4i+1})\otimes P(xx_{4i+2}) 
\xrightarrow{zero}
E(y_{8i+3})\otimes P(yy_{8i+4})
\]
\[
\xrightarrow[y_{8i+3} \ra 0]{yy_{16i+4} \ra (z_{8i+2})^2,\quad yy_{16i+12} \ra (zz_{4i+3})^4}
P(z_{8i+2})\otimes P(zz_{4i+3})
\]
 The cokernel is
$E(y_{8i+3})\otimes P(yy_{8i+4})$ and the kernel 
$E(x_{4i+1})\otimes P(xx_{4i+2}) $.
Tor is $\Gamma[w_{4i+2}]\otimes E(ww_{4i+3})$.
We have
\[
d_2(\gamma_2(w_{4i+2})) = y_{8i+3},
\]
$(w_{4i+2})^2 = yy_{8i+4}$,
and $(ww_{4i+3})^2 = w_{8i+6}$.

\item{$  \KR{3} \llra{\eta} \KR{2} \ra \K{2} \qquad $ {\bf RRK322    }}
\[
E(x_{8i+3})\otimes P(xx_{8i+4})
\xrightarrow[x_{8i+3} \ra 0]{xx_{16i+4} \ra (y_{8i+2})^2,\quad xx_{16i+12} \ra (yy_{4i+3})^4}
P(y_{8i+2})\otimes P(yy_{4i+3})
\]
\[
\xrightarrow[
(yy_{4i+3})^4 \ra 0,\quad
(y_{8i+2})^2 \ra 0]{
y_{8i+2} \ra zzz_{8i+2},\quad 
yy_{4i+3}\ra z_{4i+3}}
TP_4(z_{4i+3})\otimes E(zz_{4i}) \otimes E(zzz_{8i+2})
\]
The cokernel is $E(y_{8i+2}) \otimes TP_4(yy_{4i+3})$
and the kernel is $E(x_{8i+3})$.  Tor of this is
$\Gamma[w_{8i+4}] = E(w_{4i})$.

\item{$  \KR{2} \ra \K{2} \ra \KR{4} \qquad $ {\bf RKR224     }}
\[
P(x_{8i+2})\otimes P(xx_{4i+3})
\xrightarrow[
(xx_{4i+3})^4 \ra 0,\quad
(x_{8i+2})^2 \ra 0]{
x_{8i+2} \ra yyy_{8i+2},\quad 
xx_{4i+3}\ra y_{4i+3}}
\]\[
TP_4(y_{4i+3})\otimes E(yy_{4i}) \otimes E(yyy_{8i+2})
\xrightarrow[y_{4i+3} \ra 0,\quad yyy_{8i+2}\ra 0]{yy_{4i} \ra z_{4i}}
E(z_{4i})\otimes E(zz_{8i+5})
\]
The cokernel is $E(yy_{4i})$ and
the kernel is
$P((x_{8i+2})^2)\otimes P((xx_{4i+3})^4)$.
Tor of this is
$E(w_{16i+5}) \otimes E(ww_{16i+13})$.

\item{$  \K{2} \ra \KR{4} \llra{\eta} \KR{3}  \qquad $ {\bf KRR243     }}
\[
TP_4(x_{4i+3})\otimes E(xx_{4i}) \otimes E(xxx_{8i+2})
\xrightarrow[x_{4i+3} \ra 0,\quad xxx_{8i+2}\ra 0]{xx_{4i} \ra y_{4i}}
\]\[
E(y_{4i})\otimes E(yy_{8i+5})
\xrightarrow{zero}
E(z_{8i+3}) \otimes P(zz_{8i+4}) 
\]
The cokernel is $E(yy_{8i+5})$ and the kernel is
$TP_4(x_{4i+3})\otimes E(xxx_{8i+2})$.
Tor is 
\[
E(w_{4i}) \otimes \Gamma[ww_{16i+14}] \otimes \Gamma[www_{8i+3}]
\]
There is no $E(yy_{8i+5})$ so it
must be hit by a differential.  Rewrite Tor as
\[
E(w_{4i}) \otimes \Gamma[ww_{8i+6}] \otimes E(www_{8i+3})
\]
where the $ww_{8i+6}$ is in filtration 2.  The differential is
now obvious
\[
d_2(ww_{8i+6}) = yy_{8i+5}
\]
This leaves
\[
E(w_{4i}) \otimes E(www_{8i+3})
\]
We must have $(w_{4i})^2 =  w_{8i}$
giving us our $P(w_{8i+4})$.

\item{$   \KR{4} \llra{\eta} \KR{3} \ra \K{3}  \qquad $ {\bf RRK433      }}
\[
E(x_{4i})\otimes E(xx_{8i+5})
\xrightarrow{zero}
E(y_{8i+3}) \otimes P(yy_{8i+4}) 
\]\[
\xrightarrow[
y_{8i+3} \ra 0
]{
yy_{8i+4} \ra (zz_{4i+2})^2
}
E(z_{4i+1})\otimes P(zz_{4i+2}) 
\]
The cokernel is
$E(y_{8i+3}) \otimes P(yy_{8i+4}) $
and the kernel is
$E(x_{4i})\otimes E(xx_{8i+5})$.
Tor is $\Gamma[w_{4i+1}] \otimes \Gamma[ww_{8i+6}]$.
Rewriting the spectral sequence we have
\[
E(y_{8i+3}) \otimes P(yy_{8i+4}) \otimes 
E(w_{4i+1}) \otimes E(w_{8i+2}) \otimes 
\]\[
\Gamma[w_{16i+4}] \otimes
E(ww_{8i+6}) \otimes \Gamma[ ww_{16i+12}) 
\]
where the generating terms are in filtrations 0, 0, 1, 2, 4, 1, 2,
respectively.
To kill $y_{8i+3}$ we need two differentials
\[
d_2(ww_{16i+12}) = y_{16i+11}
\quad
\text{and}
\quad
d_4(ww_{16i+4}) = y_{16i+3}
\]
After this, what is left is
\[
 P(yy_{8i+4}) \otimes 
E(w_{4i+1}) \otimes E(w_{8i+2}) \otimes 
E(ww_{8i+6}) 
\]
Combining the last two terms we get $E(zz_{4i+2})$ with extension
$(zz_{4i+2})^2 = yy_{8i+4}$ to get our answer.

\item{$ \KR{3} \ra \K{3} \ra \KR{5}  \qquad $ {\bf RKR335       }}
\[
E(x_{8i+3}) \otimes P(xx_{8i+4}) 
\xrightarrow[
x_{8i+3} \ra 0
]{
xx_{8i+4} \ra (yy_{4i+2})^2
}
E(y_{4i+1})\otimes P(yy_{4i+2}) 
\]\[
\xrightarrow[(yy_{4i+2})^2 \ra 0]{y_{4i+1}\ra z_{4i+1},\quad yy_{4i+2}\ra zz_{4i+2} }
E(z_{4i+1})\otimes E(zz_{2i})
\]
The cokernel is
$E(y_{4i+1})\otimes E(yy_{4i+2}) $
and the kernel is 
$E(x_{8i+3}) $.
Tor is $\Gamma[w_{8i+4}]= E(w_{4i})$.

\item{$  \K{3} \ra \KR{5} \llra{\eta} \KR{4}  \qquad $ {\bf KRR354       }}
\[
E(x_{4i+1})\otimes P(xx_{4i+2}) 
\xrightarrow[(xx_{4i+2})^2 \ra 0]{x_{4i+1}\ra y_{4i+1},\quad xx_{4i+2}\ra yy_{4i+2} }
\]\[
E(y_{4i+1})\otimes E(yy_{2i})
\xrightarrow[y_{4i+1} \ra 0,\quad yy_{4i+2} \ra 0]
{yy_{4i} \ra z_{4i}}
E(z_{4i}) \otimes E(zz_{8i+5})
\]
The cokernel is
$ E(yy_{4i})$ and the kernel is
$ P((xx_{4i+2})^2) $.
Tor is $E(w_{8i+5})$.

\item{$  \KR{5} \llra{\eta} \KR{4} \ra \K{4}  \qquad $ {\bf RRK544        }}
\[
E(x_{4i+1})\otimes E(xx_{2i})
\xrightarrow[x_{4i+1} \ra 0,\quad xx_{4i+2} \ra 0]
{xx_{4i} \ra y_{4i}}
\]
\[
E(y_{4i}) \otimes E(yy_{8i+5})
\xrightarrow{zero}
TP_4(z_{4i+3})\otimes E(zz_{4i}) \otimes E(zzz_{8i+2})
\]
The cokernel is
$ E(yy_{8i+5})$ and the kernel is
$E(x_{4i+1})\otimes E(xx_{4i+2})$.
Tor is 
$
\Gamma[w_{4i+2}]
\otimes
\Gamma[ww_{4i+3}]
$.
There is no $yy_{8i+5}$ in our answer and the only differential
that could hit it is
\[
d_2(\gamma_2(ww_{4i+3})) = yy_{8i+5}
\]
This leaves us with 
$
\Gamma[w_{4i+2}]
\otimes
E(ww_{4i+3})
=
E(w_{2i})\otimes E(ww_{4i+3})
$.
We need $(ww_{4i+3})^2 =  w_{8i+6}$.

\item{$   \KR{4} \ra \K{4} \ra \KR{6}  \qquad $ {\bf RKR446         }}
\[
E(x_{4i}) \otimes E(xx_{8i+5})
\xrightarrow{zero}
TP_4(y_{4i+3})\otimes E(yy_{4i}) \otimes E(yyy_{8i+2})
\]
\[
\xrightarrow[yy_{4i} \ra (z_{2i}^2),\quad yyy_{8i+2}\ra (z_{4i+1})^2]{y_{4i+3} \ra z_{4i+3}}
TP_4(z_i)
\]
The cokernel is
$TP_4(y_{4i+3})\otimes E(yy_{4i}) \otimes E(yyy_{8i+2})$
and the kernel is
$E(x_{4i}) \otimes E(xx_{8i+5})$.
Tor is
$\Gamma[w_{4i+1}] \otimes \Gamma[ww_{8i+6}]$.
We have $(w_{4i+1})^2 = yyy_{8i+2}$ leaving, in Tor, $\Gamma[www_{4i+2}]
= E(www_{2i})$.
We have $(www_{2i})^2 = yy_{4i}$.

\item{$  \K{4} \ra \KR{6} \llra{\eta} \KR{5}  \qquad $ {\bf KRR465          }}
\[
TP_4(x_{4i+3})\otimes E(xx_{4i}) \otimes E(xxx_{8i+2})
\xrightarrow[xx_{4i} \ra (y_{2i})^2,\quad xxx_{8i+2}\ra (y_{4i+1})^2]{x_{4i+3} \ra y_{4i+3}}
\]
\[
TP_4(y_i)
\xrightarrow[y_{4i+3} \ra 0, \quad  (y_i)^2 \ra 0]
{
y_{2i} \ra zz_{2i},
\quad  y_{4i+1}\ra z_{4i+1}
}
E(z_{4i+1})\otimes E(zz_{2i})
\]
 The cokernel is
$ E(y_{2i})
\otimes
 E(y_{4i+1})$
and there is no kernel, giving us another rare short exact sequence
\[
\text{\bf S.E.S.} \qquad  \HK{4} \lra \HKR{6} \lra \HKR{5}.
\]

\item{$  \KR{6} \llra{\eta} \KR{5} \ra \K{5}  \qquad $ {\bf RRK655          }}
\[
TP_4(x_{i})
\xrightarrow[x_{4i+3} \ra 0, \: (x_i)^2 \ra 0
]
{x_{2i} \ra yy_{2i}
,  \: x_{4i+1}\ra y_{4i+1}
}
E(y_{4i+1})\otimes E(yy_{2i})
\xrightarrow { zero }
E(z_{4i+1})\otimes P(zz_{4i+2}) 
\]
There is no cokernel.  The kernel is
\[
TP_4(x_{4i+3})
\otimes
E((x_{2i})^2) 
\otimes
E((x_{4i+1})^2) 
\]
Tor, in filtration 1, is
\[
E(w_{4i})
\otimes
E(ww_{4i+1}) 
\otimes
E(www_{8i+3}) 
\]
and in filtration 2, combining all of the $\Gamma$, we get
$\Gamma[wwww_{4i+2}]$.  (This takes some manipulation.)
The $E(ww_{4i+1})$
must stay and the $E(www_{8i+3})$ must go.  To make it go,
we have a 
\[
d_3(\gamma_2(wwww_{4i+2})) = www_{8i+3}
\]
All we have left now is 
\[
E(w_{4i})
\otimes
E(ww_{4i+1}) 
\otimes
E(wwww_{4i+2}) 
\]
We have $(w_{4i})^2 = w_{8i}$
and
$(wwww_{4i+2})^2 = w_{8i+4}$.

\item{$   \KR{5} \ra \K{5} \ra \KR{7}  \qquad $ {\bf RKR557          }}
\[
E(x_{4i+1})\otimes E(xx_{2i})
\xrightarrow { zero }
E(y_{4i+1})\otimes P(yy_{4i+2}) 
\]
\[
\xrightarrow[y_{4i+1} \ra 0]{ yy_{4i+2}\ra z_{4i+2} + (zz_{2i+1})^2 }
E(z_{2i})\otimes P(zz_{2i+1})
\]
We have cokernel
$E(y_{4i+1})\otimes P(yy_{4i+2}) $
and kernel
$E(x_{4i+1})\otimes E(xx_{2i})$. 
Tor is 
$\Gamma[w_{4i+2}]\otimes \Gamma[ww_{2i+1}]$. 
We have
\[
d_2(\gamma_2(ww_{2i+1})) = y_{4i+1}
\]
and $(ww_{2i+1})^2 = yy_{4i+2}$.
The second map is unusual and doesn't come from this.
We'll pick it up in the next spectral sequence, {\bf KRR576}.

{\bf Remark}.  So little interesting happens that I should point out
when something does.  The last spectral sequence
did not completely describe the second map.  The differential is
correct, and we get the correct answer from the stated extension.
However, what really happens is 
and $(ww_{2i+1})^2 = yy_{4i+2} + w_{4i+2}$.
This doesn't change our answer, but if you look at the cokernel of
the second map, we get a $TP_4(ww_{2i+1})$ this way.  We can't see
that in the last spectral sequence, but we'll see it in the next.

\item{$   \K{5} \ra \KR{7} \llra{\eta} \KR{6}  \qquad $ {\bf KRR576          }}
\[
E(x_{4i+1})\otimes P(xx_{4i+2}) 
\xrightarrow[x_{4i+1} \ra 0]{ xx_{4i+2}\ra y_{4i+2} + (yy_{2i+1})^2 }
E(y_{2i})\otimes P(yy_{2i+1})
\]
\[
\xrightarrow[y_{4i+2} + (yy_{2i+1})^2  \ra 0]
{y_{2i} \ra (z_i)^2 ,\quad yy_{2i+1} \ra z_{2i+1}}
TP_4(z_i)
\]
We have a new problem here, and that is that the previous spectral
sequence, {\bf RKR557}, didn't pick up the details of the second map
there, i.e. the first map of this spectral sequence.  We didn't need
it there, but do here.
All we really know from the previous spectral sequence is that the
$P(xx_{4i+2})$ injects, but there are two ways to do that.  (1)
$xx_{4i+2} \ra (yy_{2i+1})^2$
or (2)
$xx_{4i+2} \ra y_{4i+2} +(yy_{2i+1})^2$.

If we try the first, our cokernel is 
$E(y_{2i})\otimes E(yy_{2i+1})$ and the kernel is $E(x_{4i+1})$.
Taking Tor of this gives $\Gamma[w_{4i+2}]$.  
There can be no $E(yy_{2i+1})$ as 
a subalgebra of $\HKR{6}$.  The only possibilities for differentials
are on $\gamma_{2^j}(w_{4i+2})$ (with $j > 0$).  
These elements are all in degrees divisible by 4 so can never hit
a $yy_{4i+1}$.  Version (1) cannot be correct, so try (2).
The cokernel now is $TP_4(yy_{2i+1})\otimes E(y_{4i})$ and the 
kernel is still $E(x_{4i+1})$ giving Tor as $\Gamma[w_{4i+2}]= w_{2i}$.  
We have $(w_{2i})^2 = y_{4i}$.

This computes the first map here and the second map in the previous
spectral sequence, {\bf RKR557}.

\item{$   \KR{7} \llra{\eta} \KR{6} \ra \K{6}  \qquad $ {\bf RRK766          }}
\[
E(x_{2i})\otimes P(xx_{2i+1})
\xrightarrow[x_{4i+2} + (xx_{2i+1})^2  \ra 0]
{x_{2i} \ra (y_i)^2 ,\quad xx_{2i+1} \ra y_{2i+1}}
TP_4(y_i)
\]
\[
\xrightarrow[y_{2i+1} \ra 0,\quad (y_i)^2 \ra 0,
\quad
y_{4i} \ra zz_{4i}
]
{
y_{8i+2}\ra zzz_{8i+2},\quad y_{8i+6} \ra (x_{4i+3})^2}
TP_4(z_{4i+3})\otimes E(zz_{4i}) \otimes E(zzz_{8i+2})
\]
The cokernel is $E(y_{2i})$ and the kernel
is 
$P( x_{4i+2} + (xx_{2i+1})^2 ) $.
Tor of this is $E(w_{4i+3})$.
We have
$(w_{4i+3})^2 = y_{8i+6}$.

\item{$   \KR{6} \ra \K{6} \ra \KR{0}  \qquad $ {\bf RKR660           }}
\[
TP_4(x_i)
\xrightarrow[x_{2i+1} \ra 0,\quad (x_{2i})^2 \ra 0,
\quad
x_{4i} \ra yy_{4i} 
]
{
x_{8i+2}\ra yyy_{8i+2},\quad x_{8i+6} \ra (x_{4i+3})^2}
TP_4(y_{4i+3})\otimes E(yy_{4i}) \otimes E(yyy_{8i+2})
\]
\[
\xrightarrow[(y_{4i+3})^2 \ra 0,\quad yyy_{8i+2}\ra 0]{y_{4i+3} \ra z_{4i+3}, \quad
yy_{4i} \ra 0 }
E(z_i) \otimes P(zz_{4i+2})
\]
The cokernel is $E(y_{4i+3})$ and the kernel is
$TP_4(x_{2i+1}) \otimes E((x_{2i})^2)$.
Tor of this is $E(w_{2i})\otimes \Gamma[ww_{8i+6}] \otimes \Gamma[www_{4i+1}]$.
We can rewrite this as algebras to be
$E(w_{2i})\otimes E(www_{4i+1}) \otimes \Gamma[ww_{4i+2}] $,
which is
$E(w_{2i})\otimes E(www_{4i+1}) \otimes E(ww_{2i}) $.
We have $(w_{2i})^2 = w_{4i}$.
This conclusion is a bit tricky.  The $ww_{2i}$ are all in higher
filtrations and for degree reasons, they could only square to
$w_{4i}$.  However, if that were the case, we would not have
enough even degree exterior generators.

\item{$  \K{6} \ra \KR{0} \llra{\eta} \KR{7}  \qquad $ {\bf KRR607           }}
\[
TP_4(x_{4i+3})\otimes E(xx_{4i}) \otimes E(xxx_{8i+2})
\xrightarrow[(x_{4i+3})^2 \ra 0,\quad xxx_{8i+2}\ra 0]{x_{4i+3} \ra y_{4i+3}, \quad
xx_{4i} \ra 0 }
\]
\[
E(y_i) \otimes P(yy_{4i+2})
\xrightarrow[y_{2i}\ra z_{2i},\quad y_{2i+1} \ra 0 ]{yy_{4i+2} \ra (zz_{2i+1})^2}
E(z_{2i})\otimes P(zz_{2i+1})
\]
The cokernel is $E(y_{4i+1})\otimes E(y_{2i}) \otimes P(yy_{4i+2})$.
The kernel is $E(xx_{4i})\otimes E(xxx_{8i+2})\otimes E((x_{4i+3})^2)$.  
Tor of this is $\Gamma[w_{2i+1}]$.
We need
\[
d_2(\gamma_2(w_{2i+1})) = y_{4i+1}
\]
To finish things off we have $(w_{2i+1})^2 = yy_{4i+2}$.

\item{$  \KR{0} \llra{\eta} \KR{7} \ra \K{7}  \qquad $ {\bf RRK077           }}
\[
E(x_i) \otimes P(xx_{4i+2})
\xrightarrow[x_{2i}\ra y_{2i},\quad x_{2i+1} \ra 0 ]{xx_{4i+2} \ra (yy_{2i+1})^2}
E(y_{2i})\otimes P(yy_{2i+1})
\]\[
\xrightarrow
[y_{2i}\ra 0, \quad yy_{4i+3} \ra 0, \quad (yy_{4i+1})^2  \ra 0]
{
yy_{4i+1} \ra z_{4i+1}
}
E(z_{4i+1})\otimes P(zz_{4i+2}) 
\]
The cokernel is $E(yy_{2i+1})$ and the kernel is
$E(x_{2i+1})$.
Tor of this is $\Gamma[w_{2i}]$.
We have
\[
d_2(\gamma_2(w_{2i})) = yy_{4i-1}
\]
This leaves $E(x_{4i+1})\otimes E(w_{2i})$.
We have $(w_{2i})^2 = w_{4i}$ to get our answer.

\item{$ \KR{7} \ra \K{7} \ra \KR{1}  \qquad $ {\bf RKR771           }}
\[
E(x_{2i})\otimes P(xx_{2i+1})
\xrightarrow
[x_{2i}\ra 0, \quad xx_{4i+3} \ra 0, \quad (xx_{4i+1})^2  \ra 0]
{
xx_{4i+1} \ra y_{4i+1}
}
\]\[
E(y_{4i+1})\otimes P(yy_{4i+2}) 
\xrightarrow[y_{4i+1} \ra 0]{ yy_{4i+2}\ra (z_{2i+1})^2 }
P(z_{2i+1})\otimes P(zz_{4i+2})
\]
The cokernel is $P(yy_{4i+2})$.
The kernel is 
$E(x_{2i}) \otimes P(xx_{4i+3}) \otimes P((xx_{4i+1})^2)$.
Tor of this is
$\Gamma[w_{2i+1}] \otimes E(ww_{4i}) \otimes E(www_{8i+3})$.
The way to untangle this mess is to have
\[
d_3(\gamma_4(w_{2i+1})) = www_{8i+3}
\]
leaving $P(yy_{4i+2})$ in filtration zero,
$E(w_{2i+1})\otimes E(ww_{4i})$
with generators in filtration 1,
and 
$E(w_{4i+2})$
with generators in filtration 2.

The only way this can work out is
\[
(w_{2i+1})^2 = yy_{4i+2}
\qquad
(ww_{4i+2})^2 = ww_{8i+4}
\qquad
(ww_{4i})^2 = ww_{8i}
\]

\item{$  \K{7} \ra \KR{1} \llra{\eta} \KR{0}  \qquad $ {\bf KRR710           }}
\[
E(x_{4i+1})\otimes P(xx_{4i+2}) 
\xrightarrow[x_{4i+1} \ra 0]{ xx_{4i+2}\ra (y_{2i+1})^2 }
P(y_{2i+1})\otimes P(yy_{4i+2})
\]
\[
\xrightarrow[(y_{2i+1})^2 \ra 0]{y_{2i+1} \ra z_{2i+1},\quad yy_{4i+2} \ra zz_{4i+2}}
E(z_{i}) \otimes P(zz_{4i+2}) 
\]
The cokernel is $E(y_{2i+1}) \otimes P(yy_{4i+2})$ and
the kernel is $E(x_{4i+1})$.  Tor of this is
$\Gamma[w_{4i+2}]=E(w_{2i})$.

\end{itemize}

\section{Appendix to the appendix: Homotopy long exact sequences}

Just for my purposes, I want to write down all the homotopy 
exact sequences associated with the fibrations from the main paper once and for
all.  
\begin{equation*}
\xymatrix@=1em{
	i 
&
\pi_i(\KU{0}) 
\ar[r]^2
&
\pi_i(\KU{0}) 
\ar[r]
&
\pi_i(\K{0})  \\
	1  & 0 & 0 & 0 \\
	0  & \Z \ar[r]^2 & \Z \ar[r] & \Zq
	}
\end{equation*}
\begin{equation*}
\xymatrix@=1em{
	i & \pi_i(\KO{1}) \ar[r]^\eta & \ar[r]  \pi_i(\KO{0}) &   \pi_i(\KU{0})  \\
	7  & 0          & 0         & 0            \\
	6  & 0          & 0         & \Z  \ar[dll]_=         \\
	5  & \Z         & 0         & 0            \\
	4  & 0          & \Z   \ar[r]^2     & \Z \ar[dll]           \\
	3  & \Zq        & 0         & 0            \\
	2  & \Zq \ar[r]^= & \Zq       & \Z \ar[dll]_2  \\
	1  & \Z \ar[r]  & \Zq       & 0            \\
	0  & 0          & \Z \ar[r]^= & \Z           
	}
\end{equation*}
\begin{equation*}
\xymatrix@=1em{
	i  & \pi_i(\KO{0}) \ar[r]^2 & \pi_i(\KO{0}) \ar[r] &  \pi_i(\KR{0}) \\
	7  & 0        & 0         & 0            \\
	6  & 0        & 0         & 0          \\
	5  & 0        & 0         & 0            \\
	4  & \Z  \ar[r]^2     & \Z  \ar[r]      & \Zq            \\
	3  & 0        & 0         & \Zq  \ar[dll]_=           \\
	2  & \Zq      & \Zq  \ar[r]^2     & \Z/(4) \ar[dll]          \\
	1  & \Zq      & \Zq  \ar[r]^=     & \Zq            \\
	0  & \Z  \ar[r]^2      & \Z \ar[r]       & \Zq           
	}
\end{equation*}
\begin{equation*}
\xymatrix@=1em{
	i  & \pi_i(\KR{1}) \ar[r]^\eta & \pi_i(\KR{0}) \ar[r] & \pi_i(\K{0}) \\
	7  & 0        & 0         & 0            \\
	6  & 0        & 0         & \Zq  \ar[dll]_=         \\
	5  & \Zq        & 0         & 0            \\
	4  & \Zq  \ar[r]^=     & \Zq        & \Zq \ar[dll]_2           \\
	3  & \Z/(4)   \ar[r]     & \Zq         &            \\
	2  & \Zq  \ar[r]^2    & \Z/(4) \ar[r]      & \Zq          \\
	1  & \Zq \ar[r]^=     & \Zq       & 0            \\
	0  & 0        & \Zq   \ar[r]^=     & \Zq           
	}
\end{equation*}

%\bibliography{/Volumes/DESKair/WORKSPACE/BIB/biblio}

\begin{thebibliography}{Moo60}

\bibitem[CS02]{CS}
D.S. Cowen and N.~Strickland.
\newblock The {Hopf} rings for ${KO}$ and ${KU}$.
\newblock {\em Journal of Pure and Applied Algebra}, 166(3):247--265, 2002.

\bibitem[EM66]{EilMoore}
Samuel Eilenberg and John~C. Moore.
\newblock Homology and fibrations. {I}. {C}oalgebras, cotensor product and its
  derived functors.
\newblock {\em Comment. Math. Helv.}, 40:199--236, 1966.

\bibitem[HK01]{HK}
P.~Hu and I.~Kriz.
\newblock Real-oriented homotopy theory and an analogue of the
  {Adams}-{Novikov} spectral sequence.
\newblock {\em Topology}, 40(2):317--399, 2001.

\bibitem[KW07]{NituER2}
N.~Kitchloo and W.~S. Wilson.
\newblock On the {Hopf} ring for {$ER(n)$}.
\newblock {\em Topology and its {A}pplications}, 154:1608--1640, 2007.

\bibitem[MM65]{MM}
J.~W. Milnor and J.~C. Moore.
\newblock On the structure of {Hopf} algebras.
\newblock {\em Annals of Mathematics}, 81(2):211--264, 1965.

\bibitem[Moo60]{Moore-bar}
J.~C. Moore.
\newblock Algèbre homologique et homologie des espaces classifiants.
\newblock In {\em Périodicité des Groupes d'Homotopie Stables des Groupes
  Classiques, d'après Bott}, volume~12 of {\em Seminaire Henri Cartan},
  chapter~7, pages 1--37. Secretariat mathematique, Paris, 1959-1960.

\bibitem[RS65]{RS}
M.~Rothenberg and N.~E. Steenrod.
\newblock The cohomology of classifying spaces of ${H}$-spaces.
\newblock {\em Bull. Amer. Math. Soc.}, 71:872--875, 1965.

\bibitem[RW77]{RW:HR}
D.~C. Ravenel and W.~S. Wilson.
\newblock The {Hopf} ring for complex cobordism.
\newblock {\em Journal of Pure and Applied Algebra}, 9:241--280, 1977.

\bibitem[Smi70]{Smith:LEMSS}
L.~Smith.
\newblock {\em Lectures on the {Eilenberg-Moore} spectral sequence}, volume 134
  of {\em Lecture Notes in Mathematics}.
\newblock Springer-Verlag, 1970.

\bibitem[Wil84]{WSW:MK}
W.~S. Wilson.
\newblock The {Hopf} ring for {Morava} {$K$}-theory.
\newblock {\em Publications of Research Institute of Mathematical Sciences,
  Kyoto University}, 20:1025--1036, 1984.

\end{thebibliography}
%\bibliographystyle{alpha}

\end{document}